\documentclass[11pt, oneside]{article} % use "amsart" for AMSLaTeX format. By default, fully justified, with COmpyter Modern font. 11pt sets the font size.

\usepackage{times} % to use Times New Roman
\usepackage[utf8]{inputenc}
\usepackage[T1]{fontenc}
\usepackage{xcolor} % to expand list of colours available: black, blue, brown, cyan, darkgrey, grey, green, lime, magenta, olive, orange, pink, purple, red, teal, violet, white, yellow
\usepackage{color} % to write in colours 
\usepackage{amsmath}
\usepackage{amssymb} % note: it loads \usepackage{amsfonts}
\usepackage{mathrsfs} % for italic letters, use the command: \mathscr{}
\usepackage{bm} % for isolated bold symbols within an equation
\usepackage{bbm} % for numbers in the \mathbb style
\usepackage{systeme} % for systems without brackets

\usepackage{amsthm} 
\theoremstyle{plain}
\newtheorem{theorem}{Theorem} % add [section] (or [chapter] in book class) to follow the section (chapter) number
\newtheorem{example}[theorem]{Example}

\usepackage{algorithm2e}
\SetKwInput{KwInputs}{Inputs}
\usepackage[noend]{algpseudocode}
\makeatletter
\def\BState{\State\hskip-\ALG@thistlm}
\makeatother

\RequirePackage[numbers]{natbib}
\RequirePackage[colorlinks,citecolor=blue,urlcolor=blue]{hyperref}
\usepackage{graphicx} % to insert images, use: \includegraphic{}
\usepackage{subcaption}

\usepackage{tikz} % to create plots with LaTeX
\usetikzlibrary{positioning ,shadows,matrix}

\graphicspath{ {./Figures/} }

\usepackage[english]{babel} % to manage the language (you can specify more than one: e.g.[english,italian], and switch among them). It automatically handles the division in syllables

\usepackage{hyperref} % for hypertext references
\hypersetup{
    pdfpagemode={UseOutlines},
    bookmarksopen,
    pdfstartview={FitH},
    colorlinks,
    linkcolor={blue}, % color of the index
    citecolor={blue}, % color of bibliographic references
    urlcolor={blue} 
}

\usepackage[tmargin=1in,bmargin=1in,lmargin=1.5in,rmargin=1.5in]{geometry} % to set the margins, use showframe to see the frame
\usepackage{setspace} % to set the spacing
\usepackage{nopageno}
\usepackage[margin=0cm]{caption}

\usepackage{enumitem} % package to set the enumerate/itemize symbols. Use, e.g., \begin{itemize}[label=$\triangleright$]
\def\iid{\overset{\textnormal{iid}}{\sim}} % i.i.d. symbol
\makeatletter  % cartesian product symbol
\@ifundefined{@currsizeindex} 
  {\RequirePackage{relsize}\let\dolarger\relsize} 
  {\def\dolarger#1{\larger[#1]}} 
\newcommand*\@@bigtimes[2]{\vphantom{\prod} 
  \vcenter{\hbox{\dolarger{4}$\m@th#1\mkern-2mu\times\mkern-2mu$}}} 
\newcommand*\bigtimes{\mathop{\mathpalette\@@bigtimes\relax}\displaylimits} 
\makeatother

\def\N{\mathbb{N}}\def\R{\mathbb{R}}\def\1{\mathbbm{1}}

\def\Acal{\mathcal{A}}\def\Fcal{\mathcal{F}}\def\Hcal{\mathcal{H}}\def\Ocal{\mathcal{O}}\def\Wcal{\mathcal{W}}

\title{\bf Bayesian Nonparametric Inference in Elliptic PDEs: Convergence Rates and Implementation}
\author{Matteo Giordano$^1$ \\ \\ $^1$University of Turin, Corso Unione Sovietica 218 bis, Turin, Italy}
\date{} % set the desired data within {}

\begin{document}

\maketitle

%\abstract{
\noindent\textbf{Abstract}\\
Parameter identification problems in partial differential equations (PDEs) consist in determining one or more functional coefficient in a PDE. In this article, the Bayesian nonparametric approach to such problems is considered. Focusing on the representative example of inferring the diffusivity function in an elliptic PDE from noisy observations of the PDE solution, the performance of Bayesian procedures based on Gaussian process priors is investigated. Building on recent developments in the literature, we derive novel asymptotic theoretical guarantees that establish posterior consistency and convergence rates for methodologically attractive Gaussian series priors based on the Dirichlet--Laplacian eigenbasis. An implementation of the associated posterior-based inference is provided and illustrated via a numerical simulation study, where excellent agreement with the theory is obtained.
%}

\bigskip

%\noindent\textbf{AMS subject classifications.} 62G20, 65N21.

%\bigskip


%\tableofcontents

%%%%%%%%%%%%%%%%%%%%%%%%%%%%%%%%%%%%%%%%%%%
\section{Introduction}
Partial differential equations (PDEs) are primary mathematical tools to model the behaviour of complex real-world phenomena, with~ubiquitous applications across engineering and the sciences. The~formulation of a PDE typically involves a number of \textit{functional parameters}, %MDPI: Please confirm if the italics should be retained. Please check the wholpe paper.
 which are often unknown in applications and not directly accessible to measurements. Employing a PDE model in practice therefore necessitates that the parameters in the equation be determined beforehand from the available data, giving rise to an \textit{inverse problem} of \textit{parameter identification}. Such problems have been extensively studied in applied mathematics~\cite{EHN96,KNS08,I17} and, more recently, in~statistics~\cite{KS04,BHM04,HP08}. See the monographs~\cite{BB18,AMOS19,N23} and the references therein for an extended overview on this research~area.

\begin{figure}%%%%%%%%%%%%%%%%%%%%%%%%%%%%%%%%%%%%
\centering
\includegraphics[width=4.5cm,height=4cm]{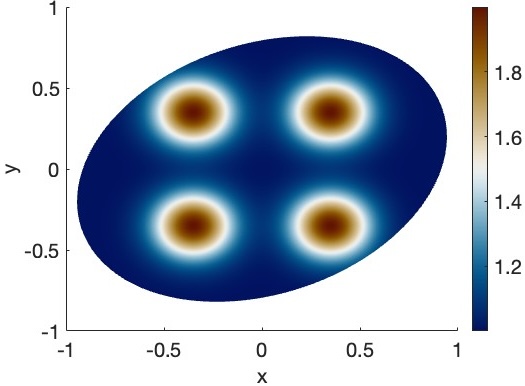}
\includegraphics[width=4.5cm,height=4cm]{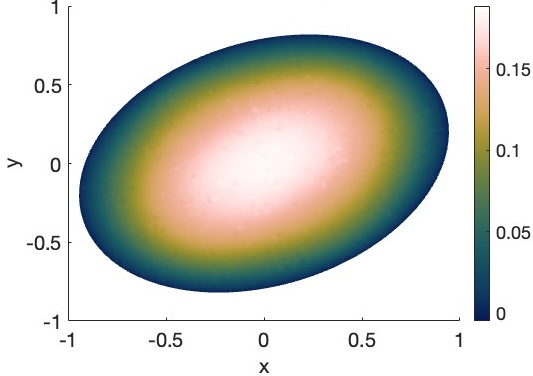}
\caption{(\textbf{Left}): An %MDPI: Please change the hyphen (-) into minus sign ($-$, "U+2212"). e.g.,~"-1" should be "$-$1".
 example of diffusivity function $f$ with four circular regions of higher conductivity. (\textbf{Right}): $n=4500$ noisy observations from the corresponding PDE solution $G(f)$.}
\label{Fig:EllipticInvProbl}
\end{figure}

	In the present paper, we shall focus on the following representative example. Consider a physical quantity undergoing \textit{diffusion} in an inhomogeneous multidimensional convex medium $\Ocal\subset\R^d, \ d\in\N$, with~smooth boundary $\partial\Ocal$. At~equilibrium, the~density $u(x)$ of the diffusing substance at any location $x\in\Ocal$ is governed by the second-order elliptic PDE
\begin{equation}
\label{Eq:EllipticPDE} 
\begin{cases}
	\nabla\cdot(f \nabla u) = s, & \text{on}\  \Ocal \\
	u=b, & \text{on}\  \partial\Ocal,
\end{cases}
\end{equation}
where $s :\Ocal\to\R$ describes the spatial distribution of local sources or sinks, $g:\partial\Ocal\to\R$ prescribes the density values at the boundary, and~the \textit{diffusivity function} $f:\Ocal\to(0,\infty)$ models spatially varying conductivity throughout the inhomogeneous domain. Under~mild regularity conditions on the PDE coefficients, standard elliptic theory implies the existence of a unique twice continuously differentiable classical solution $G(f)\equiv u_f\in C^2(\Ocal)$ to \eqref{Eq:EllipticPDE} (e.g.,~\cite{E10}, Chapter 6). Assuming that $s$ and $b$ are known, we are then interested in the problem of estimating $f$ from $n$ noisy point evaluations of $G(f)$ over a grid of (possibly random) design points $X_1,\dots,X_n$ in $\Ocal$,
\begin{equation}
\label{Eq:Obs}
	Y_i = G(f)(X_i) +  \sigma W_i, \qquad i=1,\dots,n,
\end{equation}
where $W_1,\dots,W_n$ are statistical measurement errors, and $\sigma>0$ is the noise level. In~view of the central limit theorem, the~Gaussian assumption $W_1,\dots,W_n\iid N(0,1)$ can often be realistically maintained. The~inverse problem of recovering the diffusivity in the elliptic PDE \eqref{Eq:EllipticPDE} from observations of its solution is an important building block in oil reservoir modelling~\cite{Y86}, where $u$ is the (accessible to measurements) fluid pressure, and $f$ is the (not directly observable) permeability field, which can exhibit drastic spatial variations due to changes in the reservoir composition. This problem has been studied in a large number of articles in applied mathematics, e.g.,~\cite{R81,K01,BCDPW17}, and~statistics, e.g.,~\cite{BHM04,DS11,CRSW13}. An~illustration on a bi-dimensional domain is given in Figure~\ref{Fig:EllipticInvProbl} .

%%%%%%%%%%%%%%%%%%%%%%%%%%%%%%%%%%%%%

%%%%%%%%%%%%%%%%%%%%%%%%%%%%%%%%%%%%%

%

	While the PDE \eqref{Eq:EllipticPDE} is linear, the~parameter-to-solution map $f\mapsto G(f)$ is not, which poses several methodological and theoretical challenges. In~particular, least squares functionals involving $G(f)$ are generally non-convex so that commonly used optimisation-based methods (such as Tikhonov regularisation, maximum likelihood or maximum-a posteriori estimation) cannot reliably be implemented by standard convex optimisation techniques. In~this context, the~Bayesian approach to inverse problems, popularised by influential work by Stuart~\cite{S10}, offers an attractive alternative. In~the Bayesian framework, the~unknown parameter $f$ is regarded as a random variable (with values in a function space) and assigned a \textit{prior probability distribution} $\Pi(\cdot)$ that models the available information about $f$ before collecting the observations. The~prior is then combined, through \textit{Bayes' formula}, with~the data $\{(Y_i,X_i)\}_{i=1}^n$ to form the \textit{posterior distribution} $\Pi(\cdot|\{(Y_i,X_i)\}_{i=1}^n)$, which represents the updated belief about $f$ and is used to draw the inferences. As~the posterior formally involves only evaluations of the prior and the likelihood (cf.~Equation~\eqref{Eq:Posterior} below), the approximate computation of $\Pi_n(\cdot|\{(Y_i,X_i)\}_{i=1}^n)$ and its associated \textit{posterior mean estimator} $\bar f_n := E^\Pi[f|\{(Y_i,X_i)\}_{i=1}^n]$ via sampling methods is feasible as long as the forward map $G$ can be numerically evaluated. For~ the elliptic PDE \eqref{Eq:EllipticPDE}, this can be performed using efficient PDE solvers based on finite element methods, sidestepping altogether the need for a (possibly non-existent) inversion formula for $G$, as~well as the use of optimisation approaches. In~particular, for~the class of \textit{Gaussian process priors}, efficient ad hoc Markov chain Monte Carlo (MCMC) algorithms, suited to the present infinite-dimensional setting, have been developed, e.g.,~\cite{CRSW13}. A~further decisive advantage of the Bayesian methodology is that, alongside point estimates, it also automatically delivers \textit{uncertainty quantification} for the recovery via the spread of the posterior, used in applications to provide interval-type estimators and to construct hypothesis~tests.

	The success and popularity in applications has led to recent interest in the literature for the derivation of theoretical performance guarantees for nonparametric Bayesian procedures in PDE models~\cite{V13,AN19,GN20,MNP21a,G21,AW23}. This is motivated by the fact that the performance of Bayesian methods depends on a suitable choice of the prior, which in infinite-dimensional statistical models primarily serves as a regularisation tool, and~whose specification is a delicate task in its own right (cf.~Section 1.2 in~\cite{GvdV17}). Thus, the~question arises as to whether Bayesian procedures may provide valid and prior-independent inference, at~least in the presence of informative data. The~established paradigm under which such investigation is carried out is the \textit{frequentist analysis} of Bayesian procedures, assuming that the observations are generated by a fixed \textit{ground truth} $f_0$ and studying the concentration of the posterior towards $f_0$ in the large sample size limit. We refer the reader to~\cite{GvdV17} for an introduction to this research~area.

	The present paper is concerned with the performance of Bayesian nonparametric methods based on Gaussian priors in the elliptic inverse problem \eqref{Eq:Obs}. In~Section~\ref{Sec:Methods}, we provide a general result that establishes posterior consistency and convergence rates of the conditional mean estimator for a general class of Gaussian priors (Theorem \ref{Theo:Main}). We then apply our general theorem to obtain statistical recovery rates for truncated Gaussian series priors defined on the eigenbasis of the Dirichlet--Laplacian (Example \ref{Ex:DirichLapl}), which is a commonly used basis of practical interest offering a convenient and generally applicable framework for implementation; cf.~Section \ref{Sec:DirichLapl}. This shows that the resulting procedures provide statistically valid estimation of the diffusivity $f$, with~explicit estimation error bounds that decay algebraically in the number $n$ of observations. Our results extend upon the recent investigation of Giordano and Nickl~\cite{GN20}, who considered Gaussian priors associated to popular covariance kernels (such as the Matérn or squared-exponential ones), as~well as truncated Gaussian wavelet series expansions but~did not explore the methodologically attractive procedures based on the Dirichlet--Laplacian eigenbasis constructed~herein.

\begin{figure} %%%%%%%%%%%%%%%%%%%%%%%%%%%%%%%%%%%%
\centering
\includegraphics[width=4.5cm,height=4cm]{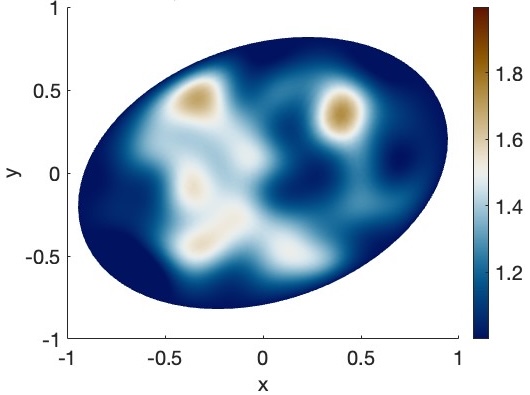}
\includegraphics[width=4.5cm,height=4cm]{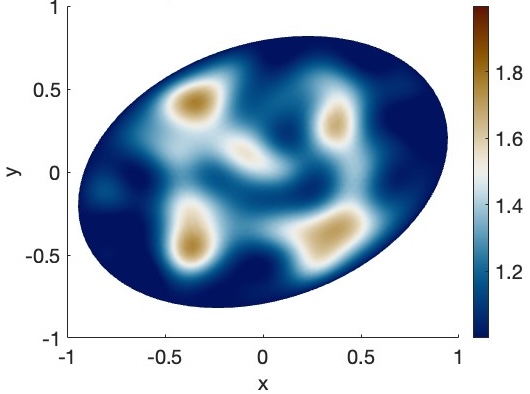}
\includegraphics[width=4.5cm,height=4cm]{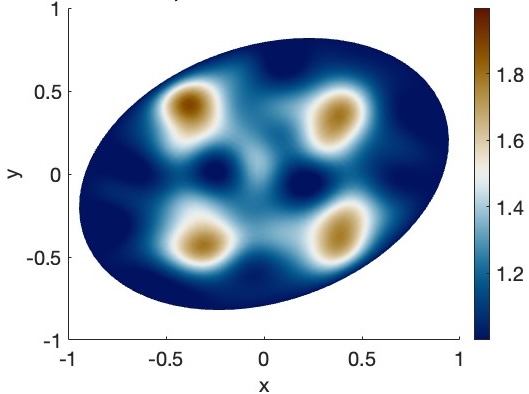}
\caption{Left %MDPI1. : Please change the hyphen (-) into minus sign ($-$, "U+2212"). e.g., "-1" should be "$-$1". 2. We moved the Figure behind its first citation. Please confirm.
 to right: posterior mean estimates $\bar f_n$ of the diffusivity function $f$ for increasing sample sizes $n= 100, 250, 1000$.}
\label{Fig:NonlinIPEstim}
\end{figure}

	Furthermore, in~Section~\ref{Sec:Results}, we complement the theoretical results with a discussion on implementation, devising two different discretisation strategies. The~first is tailored to applications where a specific set of basis functions is of interest. In~particular, in~\mbox{Section~\ref{Sec:DirichLapl},} we employ the Dirichlet--Laplacian eigenbasis, discretising the parameter space by high-dimensional truncated series expansions in accordance with our theoretical results. In~the second approach, described in Section~\ref{Sec:ImplMatern}, we discretise the parameter space via piece-wise linear functions defined on the elements of a deterministic triangular mesh, which is naturally suited to implementing Gaussian priors specified via a covariance kernel. Here, the~popular Matérn kernel is employed for illustration.  A~numerical simulation study is provided to investigate the performance of the inferential procedures under the two discretisation strategies. In~our numerical experiments, both approaches yielded satisfactory results, comparable in terms of reconstruction quality and running time. The~posterior mean estimate (relative to a Matérn process prior), computed via a Metropolis--Hastings-type MCMC algorithm, is shown in Figure~\ref{Fig:NonlinIPEstim} for increasing sample sizes, to~be compared to the true diffusion coefficient pictured in Figure~\ref{Fig:EllipticInvProbl} (left). The~reproducible MATLAB code used for the study is available at %MDPI: 1. Please confirm if the font style should be retained. 2. \hl{} %MDPI: Please add the access date (Format: Date Month Year). e.g.,~(accessed on 1 January 2020).
 \href{https://github.com/MattGiord/Bayesian-Inverse-Problems/tree/main}{\texttt{https://github.com/MattGiord}} (accessed on 19 March 2025).
%
%
%
%
%

%%%%%%%%%%%%%%%%%%%%%%%%%%%%%%%%%%%%%%%%%%
\section{Materials and~Methods}
\label{Sec:Methods}
\unskip

%%%%%%%%%%%%%%%%%%%%%%%%%%%%%%%%%%%%%%%%%%
\subsection{Likelihood, Prior and~Posterior}

Throughout, $\Ocal\subset\R^d$, $d\in\N$, is a given nonempty, open, convex and bounded set with smooth boundary $\partial\Ocal$. For~the observation model \eqref{Eq:Obs}, with~$G(f)$ as the solution to the PDE \eqref{Eq:EllipticPDE}, and~for fixed constants $\alpha>1+d/2$ and $f_{\textnormal{min}}>0$, we take the parameter space
\begin{equation}
\label{Eq:ParamSp}
\begin{split}
	&\Fcal_{\alpha,f_{\textnormal{min}}}
	 :=
		\Bigg\{f\in H^\alpha(\Ocal): 
		\textnormal{$\inf_{x \in \Ocal}f(x)\ge f_{\textnormal{min}}$, 
		$f_{|\partial \Ocal}\equiv1$, $\frac{\partial^j f}{\partial \nu^j}\equiv0$
		for $1\le j\le\alpha-1$}\Bigg\},
\end{split}
\end{equation}
with $H^\alpha(\Ocal)$ as the usual Sobolev space of regularity $\alpha$ and $\nu(x), \ x\in\partial\Ocal$, the~unit normal vector. Assuming that the source term $s$ in \eqref{Eq:EllipticPDE} is fixed and smooth, and~taking (without loss of generality) homogeneous Dirichlet boundary conditions $b\equiv 0$, the~Schauder theory for elliptic PDEs (e.g.,~Theorem 6.14 in~\cite{E10}) implies that for each $f\in \Fcal_{\alpha,f_{\textnormal{min}}}$, there exists a unique classical solution $G(f) \in C(\overline\Ocal)\cap C^{1+\alpha}(\Ocal)$ to the elliptic PDE \eqref{Eq:EllipticPDE}. We then assume data $\{(Y_i,X_i)\}_{i=1}^n$ arising as in Equation~\eqref{Eq:Obs} for some unknown $f\in \Fcal_{\alpha,f_{\textnormal{min}}}$, with~independent and identically distributed (i.i.d.) random design variables $X_1,\dots,X_n$ following the uniform distribution on $\Ocal$. Throughout, we regard the noise level $\sigma>0$ in \eqref{Eq:Obs} as fixed and known; in practice, it may often be replaced by an estimate (cf.~Section \ref{Sec:Discussion}). In~view of the i.i.d.~standard normal assumption on the noise variables $W_1,\dots,W_n$ in \eqref{Eq:Obs}, the~random vectors $\{(Y_i,X_i)\}_{i=1}^n\sim P^{(n)}_f$ have joint probability density function in product form,
$$
	p_f^{(n)}(\{(x_i,y_i)\}_{i=1}^n) 
	= \frac{1}{(2\pi\sigma^2)^{n/2}}e^{-\sum_{i=1}^n[y_i-  G(f)(x_i)]^2/(2\sigma^2)},
	\qquad y_i\in\R,\qquad x_i\in\Ocal.
$$	
Accordingly, %MDPI: Please check through the paper if indentation should be added to the first line after equations.
 the~log-likelihood is seen to be equal to, up~to an additive constant, the~negative least-square functional
\begin{equation}
\label{Eq:LogLik}
	l_n(f)
	 := -\frac{1}{2\sigma^2}\sum_{i=1}^n[Y_i-  G(f)(X_i)]^2, 
	 \qquad f\in \Fcal_{\alpha,f_{\textnormal{min}}}.
\end{equation}

	In a recent paper by Giordano and Nickl~\cite{GN20}, posterior consistency and convergence rates for the conditional mean estimator have been established for  nonparametric Bayesian procedures based on Gaussian process priors. To~incorporate the shape constraints in the parameter space $\Fcal_{\alpha,f_{\textnormal{min}}}$, ref.~\cite{GN20} employed the parametrisation
\begin{equation}
\label{Eq:Param}
	f =\Phi\circ F,
\end{equation}
where $F\in H^\alpha_0(\Ocal)$ (the completion of $C^\infty_c(\Ocal)$ with respect to $\|\cdot\|_{H^\alpha}$) and $\Phi:\R\to[f_{\textnormal{min}},\infty)$ is a \textit{regular link function}, that is a smooth, strictly increasing and bijective function with bounded derivatives and such that $\Phi(0)=1$. An~instance of regular link function is provided in Example 24 of~\cite{GN20}. In~practice, the~exponential link $\Phi(\cdot)=\exp(\cdot)$ is often used. In~the following, we occasionally switch between the notation $f$ and $F$, implicitly making use of the relation \eqref{Eq:Param}.

	Under the parametrisation \eqref{Eq:Param}, placing a prior probability distribution $\Pi(\cdot)$ on $F$ induces a (push-forward) prior on $f$, which, in~slight abuse of the notation, we still denote by $\Pi(\cdot)$. Following~\cite{GN20}, we consider \textit{scaled} Gaussian priors, constructed starting from a base (possibly $n$-dependent) centred Gaussian Borel probability measure $\Pi_{W_n}$, which we assume to be supported on a measurable linear subspace of the H\"older space $C^\beta(\Ocal)$, for~some $\beta\ge1$, and~to have reproducing kernel Hilbert space (RKHS) $\Hcal_{W_n}$ continuously embedded into $H^\alpha_0(\Ocal)$. See Chapter 2 in~\cite{GN16} for the background and terminology on Gaussian processes and measures. Given $\Pi_{W_n}$, we then construct the scaled prior $\Pi_n(\cdot)$ as the law of the random function
\begin{equation}
\label{Eq:Prior}
	V(x) := \frac{W(x)}{n^{d/(4\alpha+4+2d)}}, \qquad x\in\Ocal, 
	\qquad W\sim \Pi_{W_n}(\cdot).
\end{equation}
By linearity, $\Pi_n(\cdot)$ also defines a centred Gaussian Borel probability measure with the same support as $\Pi_{W_n}$; the scaling enforces the additional regularisation used in the theoretical analysis to deal with the nonlinearity of the inverse~problem.

	Given prior $\Pi_n(\cdot)$ as above, by~Bayes' formula (p. 7,~\cite{GvdV17}), the~posterior distribution $\Pi_n(\cdot|\{(Y_i,X_i)\}_{i=1}^n)$ of $F|\{(Y_i,X_i)\}_{i=1}^n$ arising from data in model (\ref{Eq:Obs}) equals
\begin{equation}
\label{Eq:Posterior}
	\Pi_n(B|\{(Y_i,X_i)\}_{i=1}^n)
	=
		\frac{\int_B e^{l_n(\Phi\circ F)}d\Pi_n(F)}
		{\int_{C(\Ocal)}	e^{l_n(\Phi\circ F')}d\Pi_n(F')},
	\qquad
		 B\subseteq C(\Ocal)\ \textnormal{measurable},
\end{equation}
with $l_n(\cdot)$ as the log-likelihood in \eqref{Eq:LogLik}. 

%
%
%

%%%%%%%%%%%%%%%%%%%%%%%%%%%%%%%%%%%%%%%%%%%
\subsection{Convergence~Rates}

We study the asymptotic concentration of the posterior distribution $\Pi_n(\cdot|\{(Y_i,$ $X_i)\}_{i=1}^n)$ in \eqref{Eq:Posterior} around the ground truth diffusivity function $f_0=\Phi\circ F_0$, assuming that the data $\{(Y_i,X_i)\}_{i=1}^n\sim P^{(n)}_{f_0}$ have been generated according to the observation model \eqref{Eq:Obs} with $f=f_0$. The~following theorem extends the main result of~\cite{GN20}, allowing to include in the analysis general sieve-type Gaussian priors, cf.~Example \ref{Ex:DirichLapl}. The~proof follows through similarly to Section~3.2 in~\cite{GN20}; it is included for completeness and the convenience of the reader in Appendix \ref{App:Proof}.

\begin{theorem}\label{Theo:Main}%%%%%%%%%%%%%%%%%%%%%%%%%%
	For fixed positive integers $\alpha,\beta\in\N$ such that $\alpha>\beta+d/2$, consider the scaled prior $\Pi_n$ in \eqref{Eq:Prior}, where $\Pi_{W_n}$ is a centred Gaussian Borel probability measure supported on a measurable linear subspace of $C^\beta(\Ocal)$, with~RKHS $\Hcal_{W_n}\subseteq H^\alpha_0(\Ocal)$ satisfying, for~some constant $c>0$ (independent of $n$),
$$
	\|F\|_{H^\alpha}\le c \|F\|_{\Hcal_{W_n}}, \qquad \forall F\in\Hcal_{W_n}.
$$
Further assume that 
$$
	\sup_{n\in\N}E^{\Pi_{W_n}}\|F\|_{C^1}<\infty.
$$
For fixed $F_0\in H^\alpha_0(\Ocal)$, suppose that there exists a sequence $F_{0,n}\in\Hcal_{W_n}$ such that, as~$n\to\infty$,
\begin{equation}
\label{Eq:ApproxSeq}
	\|F_0 - F_{0,n}\|_{(H^1)^*}=O(n^{-\frac{\alpha+1}{2\alpha + 2 + d}});
	\qquad \max\{\|F_{0,n}\|_{C^1},\|F_{0,n}\|_{\Hcal_{W_n}}\}=O(1).
\end{equation}
Then, there exists $L>0$ that is large enough such that, as~$n \to \infty$,
\begin{equation}
\label{Eq:PredRiskCons}
	\Pi_n \Big(f:\|G(f)-G(f_0)\|_{L^2}
	>Ln^{-\frac{\alpha+1}{2\alpha + 2 + d}}\Big|\{(Y_i,X_i)\}_{i=1}^n
	\Big) \to 0,
\end{equation}
in $P_{f_0}^{(\infty)}$-probability as $n\to\infty$. If~in addition $\beta\ge2$ and $\inf_{x \in \Ocal}s(x)$ $ >0$, then there exists $L>0$ that is large enough and a constant $\lambda>0$ such that
\begin{equation}
\label{Eq:Cons}
	\Pi_n (f:\|f-f_0\|_{L^2}
	>Ln^{-\lambda}|\{(Y_i,X_i)\}_{i=1}^n) \to 0,
\end{equation}
in $P_{f_0}^{(\infty)}$-probability as $n\to\infty$, and~moreover, the estimator $\bar f_n = \Phi\circ \bar F_n$, with~$\bar F_n = E^{\Pi_n}[F|\{(Y_i,X_i)\}_{i=1}^n]$, satisfies as $n\to\infty$,
\begin{equation}
\label{Eq:ConvRate}
	P^{(n)}_{f_0}\big(\|\bar f_n-f_0\|_{L^2}>n^{-\lambda}\big)\to0.
\end{equation}
\end{theorem}%%%%%%%%%%%%%%%%%%%%%%%%%%%%%%%%%%%%

The first statement (Equation~\eqref{Eq:PredRiskCons}) of Theorem \ref{Theo:Main} establishes posterior consistency in \textit{prediction risk}: The induced posterior on the PDE solution $G(f),\ f\sim \Pi_n (\cdot|\{(Y_i,$ $X_i)\}_{i=1}^n)$, concentrates around the true PDE solution $G(f_0)$ in $L^2$-distance at rate $n^{-(\alpha+1)/(2\alpha + 2 + d)}$. Since such a rate is known to be minimax optimal \cite{NvdGW20} (Theorem 10), procedures satisfying the assumption of Theorem \ref{Theo:Main} are seen to optimally solve the PDE-constrained regression problem of recovering $G(f_0)$ from data $\{(Y_i,X_i)\}_{i=1}^n$.

	The second statement shows that the posterior contracts around $f_0$ also in the standard $L^2$-risk, thereby solving the inverse problem of estimating the diffusivity. It follows combining \eqref{Eq:PredRiskCons} with the regularisation properties implied by the rescaling in the prior construction \eqref{Eq:Prior} and a suitable \textit{stability estimate} for $G^{-1}$. The~latter was proved in \cite{NvdGW20} (Lemma 24), and~requires the slightly stronger assumption on $\beta$ and the strict positivity of the source $s$. The~exponent $\lambda>0$ is explicitly computed in the proof of Theorem \ref{Theo:Main} and equals $\lambda=(\alpha+1)(\beta-1)/(2\alpha+2+d)(\beta+1)$. Note that $\lambda<(\alpha+1)/(2\alpha+2+d)$. While minimax optimal rates for estimating the diffusivity $f$ in model \eqref{Eq:Obs} are currently unknown, inspection of the proof shows that when $f_0\in C^\infty(\Ocal)$, then the prior can be tuned so to attain a rate as closed as desired to the parametric rate $n^{-1/2}$.

	The last statement of Theorem \ref{Theo:Main} entails that the estimator $\bar f_n$ converges towards $f_0$ in $L^2$-risk at the same rate $n^{-\lambda}$ attained by the whole posterior distribution. It is indeed a corollary of \eqref{Eq:Cons}, following from uniform integrability arguments for Gaussian measures and the Lipschitzianity of the composition with the link function $\Phi$.

%
%
%

%%%%%%%%%%%%%%%%%%%%%%%%%%%%%%%%%%%%%%%%%%
\subsection{Examples of Gaussian~Priors}

We now provide two concrete instances of Gaussian priors to which Theorem \ref{Theo:Main} applies. For~both examples, an~implementation of the associated posterior-based inference is provided in Section~\ref{Sec:Results} below. We maintain the assumption that $f_0 =\Phi\circ F_0$ for some $F_0\in H^\alpha(\Ocal)$ supported inside a given compact subset $K\subset \Ocal$. This corresponds to the common assumption that $f_0$ is known near the boundary $\partial \Ocal$ (specifically  $f_0\equiv 1$ on $\Ocal\backslash K$).

	We first consider high-dimensional Gaussian sieve priors obtained via truncating Karhunen--Loève-type random series expansions, a~frequently used approach in computation, e.g.,~\cite{DS17}. In~particular, via Theorem \ref{Theo:Main}, we study procedures based on the Dirichlet--Laplacian eigenbasis. Such constructions corresponds to commonly used Gaussian process priors with covariance kernel given by an inverse power of the Laplacian, e.g.,~\cite{S10} (Section 2.4). The~eigenbasis can be numerically computed via efficient finite element methods for elliptic eigenvalue problems, offering a broadly applicable framework for implementation on general domains $\Ocal$. Details on computation and a numerical simulation study are provided in Section~\ref{Sec:DirichLapl}.

\begin{example}[Dirichlet--Laplacian eigenbasis]\label{Ex:DirichLapl}%-----------------------------
Let $\{e_j, \ j\ge0\}\subset H^1_0(\Ocal)\cap C^\infty(\overline\Ocal)$ be the orthonormal basis of $L^2(\Ocal)$ formed by the eigenfunctions of the (negative) Dirichlet--Laplacian:
\begin{equation}
\label{Eq:Eigen}
\begin{cases}
	-\Delta e_j -\lambda_j e_j =0,  &\textnormal{on}\ \Ocal \\
	e_j=0, & \textnormal{on}\ \partial\Ocal,
\end{cases}
\qquad j\ge0,
\end{equation}
with associated eigenvalues $0<\lambda_0<\lambda_1\le\lambda_2\le\dots, $ satisfying $\lambda_j\to\infty$ as $j\to\infty$ according to Weyl's asymptotics: $\lambda_j=O(j^{2/d})$ as $j\to\infty$, cf.~Figures \ref{Fig:Eigenfuns} and \ref{Fig:Eigenvals}. See Example 6.3 and Section~7.4 in~\cite{HT07} for details. Define the Hilbert scale
$$
	\mathbb H^s
	:=\Bigg\{w\in L^2(\Ocal) \ : \ \|w\|^2_{\mathbb H^s}:=\sum_{j=0}^\infty \lambda_j^s |\langle w,e_j\rangle_2|^2
	<\infty\Bigg\},
	\qquad s\ge0.
$$

We then have $\mathbb H^0=L^2(\Ocal)$ (with equality of norms), $\mathbb H^1= H^1_0(\Ocal)$ (with equivalence of norms), $\mathbb H^2= H^1_0(\Ocal)\cap H^2(\Ocal)$, and~the continuous embedding $\mathbb H^s\subset H^s(\Ocal)$ for all $s\ge0$ (holding generally with strict inclusion), cf.~p.~472f.~in~\cite{T11}. In~fact, it holds that $\|w\|_{\mathbb H^s}\simeq \|w\|_{H^s}$ for all $w\in \mathbb H^s$ and $s\ge0$ (proved initially for $s=2$, then extended by induction to all larger integers, and~finally by interpolation to all positive reals), and~if $F\in H^s(\Ocal)$ is compactly supported within $\Ocal$, then $F\in \mathbb H^s$. Finally, defining $\mathbb H^{-1}:=(\mathbb H^1)^*=(H^1_0(\Ocal))^*$, we have the equivalence (cf.~Equation~(A15) in~\cite{T11})
\begin{align}
\label{Eq:DualNorm}
	\|w\|^2_{(H^1_0)^*}
	\simeq \|w\|^2_{\mathbb H^{-1}}
	:= \sum_{j=0}^\infty \lambda_j^{-1} |\langle w,e_j\rangle_2|^2.
\end{align}

\vspace{-6pt}
\begin{figure} %%%%%%%%%%%%%%%%%%%%%%%%%%%%%%%%%%%
\centering
\includegraphics[width=4.5cm,height=4cm]{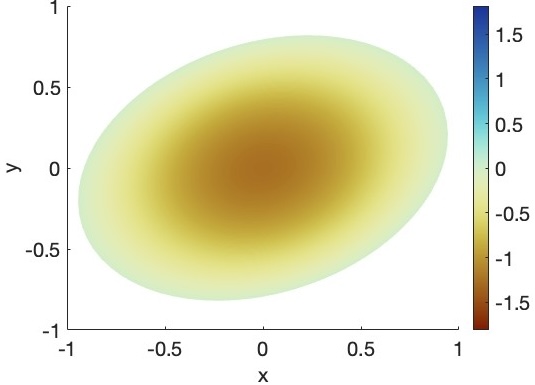}
\includegraphics[width=4.5cm,height=4cm]{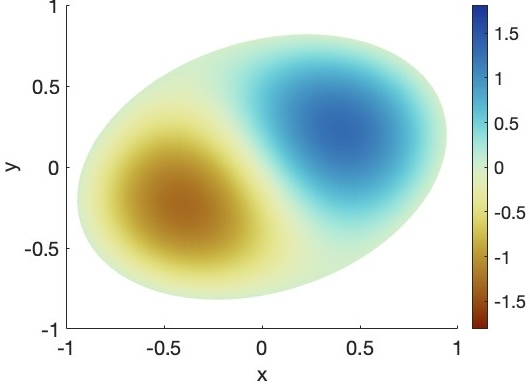}
\includegraphics[width=4.5cm,height=4cm]{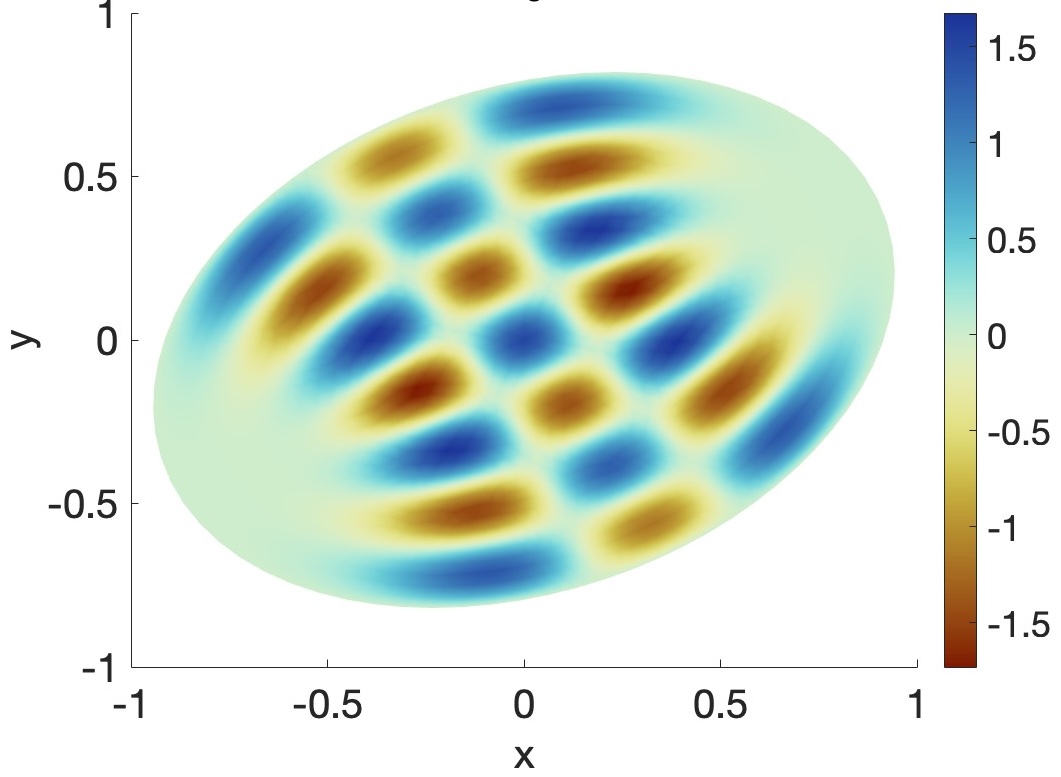}
\caption{Left %MDPI: Please change the hyphen (-) into minus sign ($-$, "U+2212"). e.g., "-1" should be "$-$1".
 to right: The first, second, and fiftieth Dirichlet--Laplacian eigenfunctions $e_0, \ e_1$ and $e_{49}$, computed via finite element~methods.}
\label{Fig:Eigenfuns}
\end{figure}%%%%%%%%%%%%%%%%%%%%%%%%%%%%%%%%%%%
\vspace{-6pt}

\begin{figure}[t] %%%%%%%%%%%%%%%%%%%%%%%%%%%%%%%%%%%
\centering
\includegraphics[width=4.5cm,height=4cm]{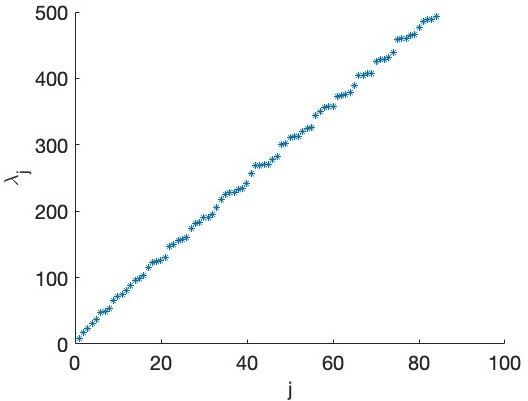}
\caption{Dirichlet--Laplacian eigenvalues $\lambda_j$ in the range $[0,500]$, computed via finite element~methods.}
\label{Fig:Eigenvals}
\end{figure}%%%%%%%%%%%%%%%%%%%%%%%%%%%%%%%%%%%

	Now for fixed $J\in\N$, the~Gaussian random sum
$$
		\overline W_J:= \sum_{j\le J}\lambda_j^{-\alpha/2} W_j e_j,
		 \qquad W_j\iid N(0,1),
$$
defines a centred Gaussian Borel probability measure supported on (and with RKHS equal to) the finite-dimensional space $\Hcal_{\overline W_J}:=\textnormal{span}\{ e_j, \ j\le J \}
%	\subset H^1_0(\Ocal)\cap C^\infty(\overline\Ocal)\cap \bigcap_{s=0}^\infty\mathbb H^s,
$, with~RKHS norm
$$
	\|\bar w\|^2_{\Hcal_{\overline W_J}}=\sum_{j\le J}\lambda_j^{\alpha} w_j^2
	= \|\bar w\|^2_{\mathbb H^\alpha}\simeq\|\bar w\|^2_{H^\alpha}, \qquad \forall w\in \overline W_J.
$$

	Fix any smooth cut-off function $\chi\in C^\infty_c(\Ocal)$ such that $\chi=1$ on $K$, and~consider the random function
\begin{equation}
\label{Eq:DirichLaplBasePrior}
		 W_n := \chi W_{J_n}= \sum_{j\le J_n}\lambda_j^{-\alpha/2} W_j \chi e_j,
		 \qquad W_j\iid N(0,1),
\end{equation}
where $J_n \in \mathbb N$ is such that $J_n \simeq n^{d/(2\alpha +2 +d)}$. By~the linearity and boundedness of multiplication by $\chi$, the~law $\Pi_{W_n}$ of $W_n$ defines, according to Exercise 2.6.5 in~\cite{GN16}, a~centred Gaussian prior supported on (and with RKHS equal to)
$$
	\Hcal_{W_n}=\textnormal{span}\{\chi e_j, \ j\le J_n \}
	\subset C^\infty_c(\Ocal) \subset
	\bigcap_{s=0}^\infty H^s_0(\Ocal)
	\cap\bigcap_{s=0}^\infty \mathbb H^s,
$$
with RKHS norm satisfying, with~a multiplicative constant independent of $n$,
$$
	\| \chi \bar w\|_{\Hcal_{W_n}}\le \|\bar w\|_{\Hcal_{\overline W_{J_n}}} 
	\simeq \|\bar w\|_{H^\alpha},
	\qquad \forall \bar w\in \Hcal_{\overline W_{J_n}}.
$$
Arguing as in Example 25 in~\cite{GN20}, one further shows that for some constant $c>0$ (independent of~$n$),
$$
	\|w\|_{H^\alpha}\le c \|w\|_{\Hcal_{W_n}}, \qquad
	\forall  w\in \Hcal_{W_n}.
$$

	Finally, by~a Sobolev embedding and the above inequality,
\begin{align*}
	E^{\Pi_{W_n}}\|F\|^2_{C^1}
	&\lesssim E^{\Pi_{W_n}}\|F\|^2_{H^\alpha}
	\le c E^{\Pi_{W_n}}\|F\|^2_{\Hcal_{W_n}}
	\le E\left[ \sum_{j\le J_n} \lambda_j^{-\alpha}W_j^2 \right],
\end{align*}
which is uniformly bounded in $n$ recalling that $W_i\iid N(0,1)$ and the fact that by Weyl's asymptotics  $\lambda_j^{-\alpha}=O(j^{-2\alpha/d})$ with $\alpha>1+d/2$. This shows that the sequence of base priors $\Pi_{W_n}$ satisfies the first two assumptions of Theorem \ref{Theo:Main}. For~ground truths $F_0\in H^\alpha(\Ocal)$ compactly supported inside $K$, we have $F_0\in\mathbb H^\alpha$. Construct the finite-dimensional approximations
\begin{equation}
\label{Eq:Approx}
	F_{0,n}=\sum_{j\le J_n} \langle F_0,e_j\rangle_2\chi e_j \in \Hcal_{W_n}, \qquad n\in\N.
\end{equation}
Then for all $n\in\N$,
\begin{align*}
	\|F_{0,n} \|_{\Hcal_{W_n}}
	&\le \left\|\sum_{j\le J_n} \langle F_0,e_j\rangle_2 e_j
	\right\|_{\Hcal_{\overline W_{J_n}}}
	=\left\|\sum_{j\le J_n} \langle F_0,e_j\rangle_2 e_j
	\right\|_{\mathbb H^\alpha}
	\le \|F_0\|_{\mathbb H^\alpha}
	\simeq  \|F_0\|_{H^\alpha}<\infty.
\end{align*}
By a Sobolev embedding, we similarly have $\|F_{0,n} \|_{C^1}\le  \|F_0\|_{H^\alpha}<\infty$ for all $n\in\N$. Furthermore, since both $F_0$ and $F_{0,n}$ have compact support within $\Ocal$,
\begin{align*}
	\|F_0 - F_{0,n}\|_{(H^1)^*}
	&=\sup_{H\in H^1(\Ocal)}\int_\Ocal(F_0(x) - F_{0,n}(x))H(x)dx\\
	&=\sup_{H\in H^1_0(\Ocal)}\int_\Ocal(F_0(x) - F_{0,n}(x))H(x)dx
	=\|F_0 - F_{0,n}\|_{(H^1_0)^*},
\end{align*}
and recalling \eqref{Eq:DualNorm}, Weyl's asymptotics, and~the choice of $J_n$,
\begin{align*}
	\|F_0 - F_{0,n}\|_{(H^1)^*}^2
	&=\sum_{j>J_n}\lambda_j^{-(1+\alpha)}
	\lambda_j^{\alpha}|\langle F_0,e_j\rangle_2|^2\\
	&\le \lambda_{J_n}^{-(1+\alpha)}\|F_0\|_{\mathbb H^\alpha}^2
	\lesssim (J_n)^{-2(1+\alpha)/d}\simeq n^{-2(\alpha+1)/(2\alpha+2+d)}.
\end{align*}
We conclude that Theorem \ref{Theo:Main} applies with the sequence of base Gaussian sieve priors in \eqref{Eq:DirichLaplBasePrior}, choosing the approximations $F_{0,n}$ according to \eqref{Eq:Approx}.
\end{example}%%%%%%%%%%%%%%%%%%%%%%%%%%%%%%%%%%%%

	The second main example, already considered in~\cite{GN20}, concerns stationary Gaussian processes specified via a translation invariant covariance kernel. For~concreteness, we focus on the popular Matérn kernel. Implementation of the resulting procedures is illustrated in Section~\ref{Sec:ImplMatern}.

\begin{example}[Matérn covariance kernel]\label{Ex:Matern}%%%%%%%%%%%%%%%%%%%%%%%%%%%%%%%%%%%
	Consider a Mat\'ern process $W=\{W(x),\ x\in\Ocal\}$ on $\Ocal$ with regularity $\alpha-d/2$, that is a centred stationary Gaussian process with covariance kernel
\begin{equation}
\label{Eq:MatCov}
	C(x,y)= \frac{2^{1-\alpha}}{\Gamma(\alpha)}\left(\frac{|x-y|\sqrt {2 \alpha}}{\ell}\right)^\alpha
	B_\alpha \left(\frac{|x-y|\sqrt{2\alpha}}{\ell}\right),
	\qquad
	x,y\in\Ocal,
	\qquad 
	\ell>0,
\end{equation}
where $\Gamma$ denotes the gamma function and $B_\alpha$ is the modified Bessel function of the second kind. Fix any smooth cut-off function $\chi\in C^\infty_c(\Ocal)$ such that $\chi=1$ on $K$. It can then be shown (cf.~Example 25 in~\cite{GN20})  that the law $\Pi_W$ of $\chi W$ defines a centred Gaussian Borel probability measure supported on the separable linear subspace $C^{\beta'}(\Ocal)$ of $C^\beta(\Ocal)$ for any $\beta<\beta'<\alpha-d/2$. Furthermore, its RKHS is given by $\Hcal_W=\{\chi F,\ F\in H^\alpha(\Ocal)\}\subset H^\alpha_0(\Ocal)$, with~the RKHS norm satisfying
$$
	\|\chi F\|_{H^\alpha}\lesssim \|\chi F\|_{\Hcal_W}, \qquad \forall F\in H^\alpha(\Ocal).
$$
For ground truths $F_0\in H^\alpha(\Ocal)$ compactly supported inside $K$, we have $\chi F = F$, so that $F\in \Hcal_W$. We conclude that Theorem \ref{Theo:Main} applies for a base Matérn process prior $\Pi_{W_n}:=\Pi_{W}$, choosing the trivial approximating sequence $F_{0,n}:=F_0$ for all $n\in\N$. 
\end{example}%%%%%%%%%%%%%%%%%%%%%%%%%%%%%%%%%%%%

%
%
%
%
%

%%%%%%%%%%%%%%%%%%%%%%%%%%%%%%%%%%%%%%%%%%
\section{Results}\label{Sec:Results}

We investigate the performance of the procedures based on the Gaussian priors considered in Examples \ref{Ex:DirichLapl} and \ref{Ex:Matern} in a numerical simulation study. For~illustration, we fix the working domain $\Ocal$ to the area contained inside the rotated ellipse, with~$\theta=\pi/6$,
$$
	\partial\Ocal=
	\{(\cos(t)\cos(\theta) - 0.75\sin(t)\sin(\theta),0.75\sin(t)\cos(\theta)+\cos(t)\sin(\theta)),\
	t\in[0,2\pi)\},
$$
and take the ground truth conductivity
\begin{equation}
\label{Eq:f0}
\begin{split}
	f_0(x,y)
	&= \sum_{k,m\in\{-1,1\}}e^{-(5x+1.75k)^2-(5y+1.75m)^2},
	\qquad (x,y)\in\Ocal,
\end{split}
\end{equation}
cf.~Figure \ref{Fig:EllipticInvProbl} (left). We then generate observations $\{(Y_i,X_i)\}_{i=1}^n$ according to \eqref{Eq:Obs}. The~true PDE solution $G(f_0)$ is numerically computed via finite element methods, using the MATLAB PDE Toolbox. For~the experiments, the~source function is set to  $s(x,y) = \exp(-(5x-2.5)^2-(5y)^2)+\exp(-(7.5x)^2-(2.5y)^2)+\exp(-(5x+2.5)^2-(5y)^2), \ (x,y)\in\Ocal$, and~the noise standard deviation to $\sigma=0.0025$ (with corresponding signal-to-noise ratio \mbox{$\|G(f_0)\|_{L^2}/\sigma =68.50$}). Further experiments, with~differently shaped ground truths, are presented in Appendix \ref{App:MoreNum}.

%
%
%

%%%%%%%%%%%%%%%%%%%%%%%%%%%%%%%%%%%%%%%%%%
\subsection{Results with Truncated Gaussian Series~Priors}
\label{Sec:DirichLapl}

To implement the truncated Gaussian series priors from Example \ref{Ex:DirichLapl}, we numerically compute the first $J_n\simeq n^{d/(2\alpha+2+d)}$ Dirichlet--Laplacian eigenpairs via finite element methods; see Figures~\ref{Fig:Eigenfuns} and \ref{Fig:Eigenvals}. Identifying  the functional parameter $F$ in \eqref{Eq:Param} with the vector of coefficients %MDPI: Please check if the bold in the equations are necessary and keep consistent in the whole text.
 $\mathbf{F} := [F_1,\dots,F_{J_n}]$, where $F_j:=\langle F,e_j\rangle_2$, the~prior then corresponds to
\begin{equation}
\label{Eq:DiscrPrior}
	\mathbf{F}\sim N(0,\Lambda_n), \qquad \Lambda_n
	:=n^{-d/(2\alpha+2+d)}
	\textnormal{diag}(\lambda_1^{-\alpha},\dots,\lambda_{J_n}^{-\alpha})\in\R^{J_n,J_n}.
\end{equation}

	While the employed prior is Gaussian, the~nonlinearity of the forward map $f\mapsto G(f)$ implies that the log-likelihood $l_n(f)$ in \eqref{Eq:LogLik} depends nonlinearly on $f$. The~resulting posterior distribution is therefore generally non-Gaussian and not available in closed form. We then resort to MCMC sampling via the pCN algorithm~\cite{CRSW13}, which is a specific instance of the random-walk Metropolis--Hastings method for Gaussian priors that is known to be robust to the discretisation dimension. In~the present setting, we employ the pCN algorithm to generate an $\R^{J_n}$-valued Markov chain $(\mathbf{F}_h)_{h\ge1}$ with invariant measure equal to the posterior distribution, starting from an initialisation point $\mathbf{F}_0$ and then, for~$h=0,1,2,\dots$, repeating the following~steps:
\begin{enumerate}
	\item Draw a prior sample $\xi\sim N(0,\Lambda_n)$, where $\Lambda_n$ is as in \eqref{Eq:DiscrPrior}, and~for $\delta>0$ define the proposal $p := \sqrt{1-2\delta}\mathbf{F}_h + \sqrt{2\delta}\xi$;

\item Set
$$
	\mathbf{F}_{h+1} :=
	\begin{cases}
	p, & \textnormal{with probability}\ 1\land
	e^{l(\Phi\circ p)-l(\Phi\circ \mathbf{F}_h)},\\
	\mathbf{F}_h, & \textnormal{otherwise},
	\end{cases}
$$
where $l$ is the log-likelihood function in \eqref{Eq:LogLik}.
\end{enumerate}
For each iteration, step 2~requires the evaluation of the log-likelihood $l_n(\Phi\circ p)$, which in turn entails the numerical evaluation of the PDE solution $G(\Phi\circ p)$ at the design points $X_1,\dots,X_n$, cf.~\eqref{Eq:LogLik}, which we perform via finite element methods. The~prior samples $\xi\sim N(0,\Lambda_n)$ required in step 1 are straightforward to draw since $\Lambda_n$ is~diagonal.

	The algorithm is terminated after $H$ steps, returning approximate samples $\{\mathbf{F}_0,\ h=0,\dots, $ $H\}$ from the posterior distribution, where a first batch of iterates is typically discarded as the burn-in. Under~a set of assumptions on the forward map that is compatible with the present setting, Hairer, Stuart, and Vollmer~\cite{HSV14} derived dimension-free spectral gaps which imply rapid convergence of the pCN marginal laws towards the invariant measure. As~a consequence, the~posterior mean $\bar F_n$ can reliably be numerically computed by the MCMC~average
\begin{equation}
\label{Eq:MCMCAverage}
	\bar{\mathbf{F}} := \frac{1}{H+1}\sum_{h=0}^H \mathbf{F}_h,
\end{equation}
with non-asymptotic bounds for the numerical approximation error. Posterior credible sets can likewise be reliably computed by considering the empirical quantiles of the pCN~samples.

\begin{figure} %%%%%%%%%%%%%%%%%%%%%%%%%%%%%%%%%%%%%%%%
\centering
\includegraphics[width=4.5cm,height=4cm]{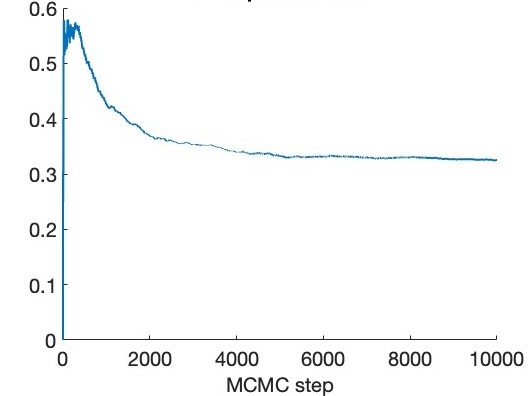}
\includegraphics[width=4.5cm,height=4cm]{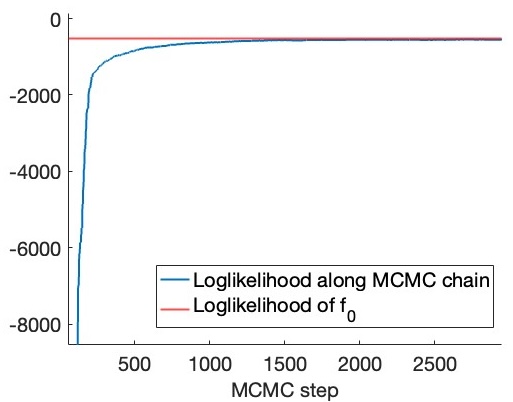}
\caption{(\textbf{Left}): %MDPI: Please change the hyphen (-) into minus sign ($-$, "U+2212"). e.g., "-1" should be "$-$1". 2. Please use commas to separate thousands for numbers with five or more digits (not for four digits) in the picture. e.g., "10000" should be "10,000"
 acceptance rate over the first 10,000 pCN samples. The~rate stabilises around $30\%$ after the initial burn-in time (first 5000 iterates). (\textbf{Right}): in blue, the~log-likelihood $l_n(\Phi\circ\vartheta_h)$ of the first 3000 iterates; in red, the~log-likelihood $l_n(f_0)$ of the true diffusion coefficient $f_0$.}
\label{Fig:NonlinIPAcceptLoglik}
\end{figure}

 Figure~\ref{Fig:NonlinIPEstim} shows the posterior mean estimates $\bar f_n = \Phi\circ \bar F_n$ of the conductivity function $f$ obtained for increasing sample sizes $n=100,250,1000$. The~corresponding $L^2$ estimation errors (and the associated relative errors) are reported in Table~\ref{Tab:Errs} (first and second row), displaying a progressive decay as expected from the theory developed in Section~\ref{Sec:Methods}. Across the experiments, the~prior regularity parameter $\alpha$ in Equation~\eqref{Eq:DiscrPrior} is set to \mbox{$\alpha=3=2+d/2$.} 
The pCN algorithm is iterated $H$ = 25,000 times, and~the first $5000$ samples are discarded as the burn-in. The~step size is set to %MDPI: We added "0" before the decimal sign. Please confirm this revision.
$\delta =  0.0025, 0.001, 0.0005, 0.00025$, tuned depending on the sample size to obtain a stabilisation of the acceptance rate around 30\% after the burn-in; see Figure~\ref{Fig:NonlinIPAcceptLoglik} (left). A~pCN \textit{cold-start} $\mathbf{F}_0 = 0$ is employed. Figure~\ref{Fig:NonlinIPAcceptLoglik} (right) shows, for~the experiment with $n=1000$, the~log-likelihood $l_n(\Phi\circ\vartheta_h)$ along the first $3000$ pCN iterates (part of the burn-in phase), seen to rapidly increase towards, and~then stabilise around, the~log-likelihood $l_n(f_0)$ attained by the true diffusivity $f_0$. The~computation time per experiment is around 50 min on a 2020 M1 MacBook Pro 13 running MATLAB~2024b.

\begin{table} 
\centering
\caption{$L^2$-estimation %MDPI: 1. We removed the vertical lines. Please confirm. 2. Please check if the table has no table header.
 errors and relative errors for the posterior mean estimates with increasing sample~sizes.}%%%%%%%%%%%%%%%%%%%%%%%%%%%%%%%%%%%%%%%%
%\centering
\begin{tabular}{ccccc} 
\hline
 	$n$ 			 &   100  & 250   & 500 & 1000 \\
\hline
 $\|\bar f_n - f_0\|_{L^2}$; series prior &  0.2981 &    0.2232  &  0.2144 & 0.1581\\
\hline
 $\|\bar f_n - f_0\|_2/\|f_0\|_{L^2}$; series prior &  17.67\% &  12.23\% &   12.71\%  & 9.36\%\\
\hline
 $\|\bar f_n - f_0\|_{L^2}$; Matérn prior & 0.3289  &    0.2677  &  0.2033  & 0.1647\\
\hline
 $\|\bar f_n - f_0\|_2/\|f_0\|_{L^2}$; Matérn prior &  18.98\% &  15.86\% &  12.05\%   & 9.76\%\\
\hline
\end{tabular}
\label{Tab:Errs}
\end{table}%%%%%%%%%%%%%%%%%%%%%%%%%%%%%%%%%%%%%%%%%

\vspace{-6pt}
%%%%%%%%%%%%%%%%%%%%%%%%%%%%%%%%%%%%%%%%

%
%
%

%%%%%%%%%%%%%%%%%%%%%%%%%%%%%%%%%%%%%%%%%%
\subsection{Results with the Matérn Process~Prior}
\label{Sec:ImplMatern}

For Gaussian process priors specified via a covariance kernel, such as the Matérn process considered in Example \ref{Ex:Matern}, we discretise the parameter space by assuming that $F$ in \eqref{Eq:Param} is given by the finite sum
\begin{equation}
\label{Eq:Sum2}
	F = \sum_{m=1}^M F_{m}\psi_m, \qquad F_{m}\in\R,
	\qquad M\in\N,
\end{equation}
where $\{\psi_1 ,\dots, \psi_M\}$ are piecewise linear functions on a deterministic triangular mesh with nodes $\{z_1,\dots,z_M\}\subset\Ocal$, uniquely characterised by the relation $\psi_m(z_{m'}) = 1_{\{m=m'\}}$. Accordingly, $F$ in \eqref{Eq:Sum2} satisfies $F(z_m) = F_m$, and~for any $x\in\Ocal$, the value $F(x)$ is obtained by linear interpolation over the pairs $\{(z_m,F_m), \ m=1,\dots,M\}$. Given the Matérn kernel $C$ in \eqref{Eq:MatCov}, and~identifying $F$ with the vector of values $\mathbf{F}:=[F_1,\dots,F_M]$, the~prior then corresponds to
$$
	\mathbf{F}\sim N(0,\mathbf{C}),
	\qquad
	\mathbf{C} = [\textnormal{C}_{mm'}]_{m,m'=1}^M\in\R^{M,M}, 
	\qquad \textnormal{C}_{mm'}=C(z_m,z_m').
$$

	Posterior inference based on the Matérn process prior may be implemented via the pCN algorithm described in Section~\ref{Sec:DirichLapl} above. For~each iteration, the~construction of the pCN proposal $p$ in step 1~entails sampling an independent multivariate Gaussian random variable $\xi\sim N(0,\mathbf C)$. In~step 2, the~computation of the acceptance probability requires the evaluation of the proposal log-likelihood $l_n(\Phi\circ p)$, which can be performed via numerical PDE methods as described~above.

	For the ground truth $f_0$ specified in Equation~\eqref{Eq:f0}, Table~\ref{Tab:Errs} (third and fourth row) displays the $L^2$-estimation errors (and the relative errors) associated to the posterior mean estimates $\bar f_n = \Phi\circ \bar F_n$ for increasing sample sizes $n=100,250,500,1000$. Across the experiments, the~parameter space is discretised using a uniform triangular mesh with $M = 2000$ nodes. The~hyperparameters for the Matérn covariance kernel in \eqref{Eq:MatCov} are set to $\alpha = 3$ and $\ell = 0.25$. Similarly to the results in Section~\ref{Sec:DirichLapl}, the~runs of pCN are stopped after \mbox{$H$ = 25,000}~iterations, with~a tuning of the step size (respectively, $\delta =  0.0025, 0.001, 0.0005, 0.00025$) to achieve a stabilisation of the acceptance rate at around $30\%$ after the burn-in (corresponding to the first 5000 iterates). A~non-informative initialisation point $\mathbf F_0 = 0$ is chosen for each run. The~computation times ranges between 50 and 57 min, in~line with those obtained for the truncated Gaussian series prior in Section~\ref{Sec:DirichLapl}.

%
%
%
%
%

%%%%%%%%%%%%%%%%%%%%%%%%%%%%%%%%%%%%%%%%%%%%%%%
\section{Discussion}\label{Sec:Discussion}

In this article, we have considered the nonparametric Bayesian approach to inference in elliptic PDEs, focusing on the standard benchmark problem of estimating the diffusivity function from noisy observation of the PDE solution. We have provided a general asymptotic concentration result, Theorem \ref{Theo:Main}, for~the posterior distribution and the posterior mean estimator, and~showed that it applies to two classes of Gaussian priors of interest, namely truncated Gaussian series priors defined on the Dirichlet--Laplacian eigenbasis (cf.~Example~\ref{Ex:DirichLapl}) and Matérn process priors (cf.~Example \ref{Ex:Matern}). For~both prior models, we have devised strategies for implementing the posterior-based inference, employing efficient and reliable MCMC algorithms.  The~performance of the considered methods has been investigated in a numerical simulation study, where excellent reconstruction results and a general agreement with the theory have been obtained. Overall, the~main advantages of the approach lie in its modelling flexibility, and the~robustness of the associated computational methods, as~well as the availability of theoretical guarantees on the recovery~performance.

	Let us conclude overviewing some related research questions. Firstly, we remark that that the noise standard deviation $\sigma$ in \eqref{Eq:Obs} is assumed to be known throughout the paper. In~the realistic scenario where $\sigma$ is also unknown, the~methodology developed here can readily be adapted by replacing it (in an \textit{empirical Bayes} spirit) with a preliminary (non-likelihood-based) estimate. Several strategies have been proposed in the literature for variance estimation in nonparametric regression models, ranging from residual-based estimators using kernel smoothing~\cite{hall1990variance} and splines~\cite{wahba1978improper}, to~difference-based estimators~\cite{rice1984bandwidth}. Alternatively, a~joint Bayesian model for $f$ and $\sigma$ in \eqref{Eq:Obs} could be considered by endowing $\sigma$ with a prior distribution (for example, an~independent conjugate inverse gamma prior). This leads to the interesting question concerning the extension of the theoretical results presented in Section~\ref{Sec:Methods} to the setting with unknown variance; see~\cite{kejzlar2021fast} for related results in a direct regression~model.

	Secondly, we mention the important issue of specifying the hyperparameter values for the considered prior distributions, namely, the regularity parameter $\alpha$ for the truncated Gaussian series prior from Example \ref{Ex:DirichLapl}, and~the smoothness and length-scale in the Matérn covariance kernel \eqref{Eq:MatCov}. There a vast literature investigating the methodological and theoretical aspects of empirical and hierarchical Bayesian strategies to fully data-driven hyperparameters selection; see~\cite{KSvdVvZ15, roussszabo,teckentrup2020convergence,agapiou2014analysis,GN20}, where many more references can be found. Investigating the implications and performance of these methods in the context of the elliptic PDE model and the Gaussian prior distributions considered in the present article is an interesting direction for future~research. 
	
\appendix

%%%%%%%%%%%%%%%%%%%%%%%%%%%%%%%%%%%%%%%%%%
\section[\appendixname~\thesection]{Further Numerical Results}
\label{App:MoreNum}

\begin{figure}[t]
\centering
\includegraphics[width=4.5cm,height=4cm]{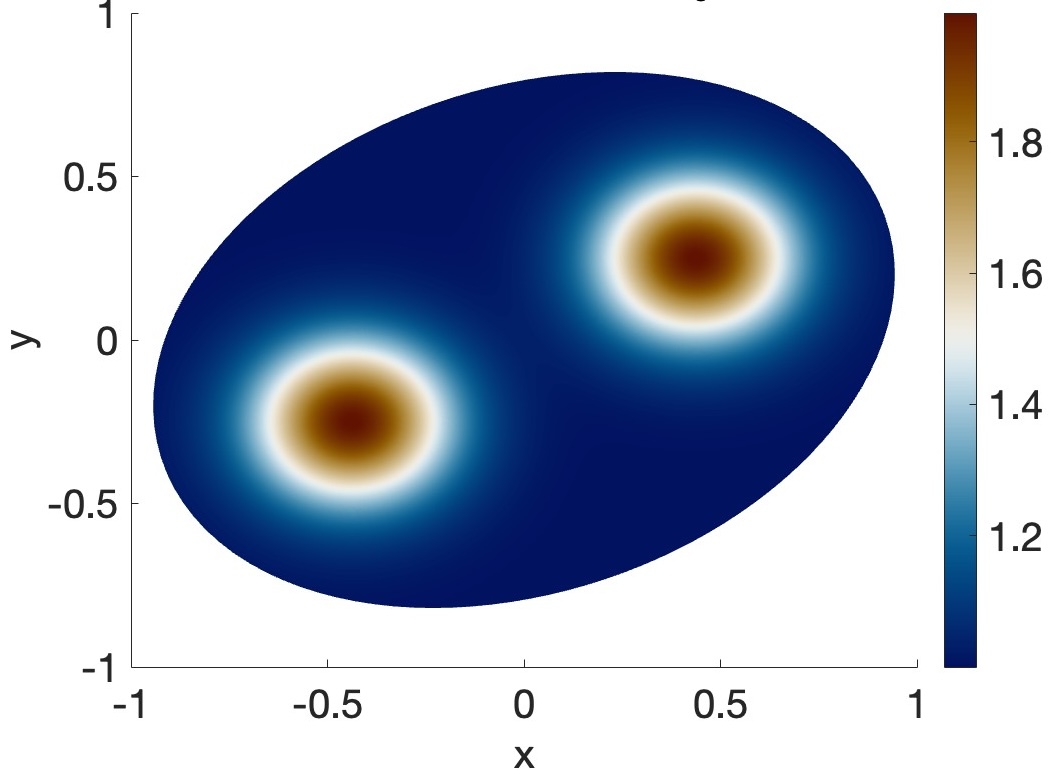}
\includegraphics[width=4.5cm,height=4cm]{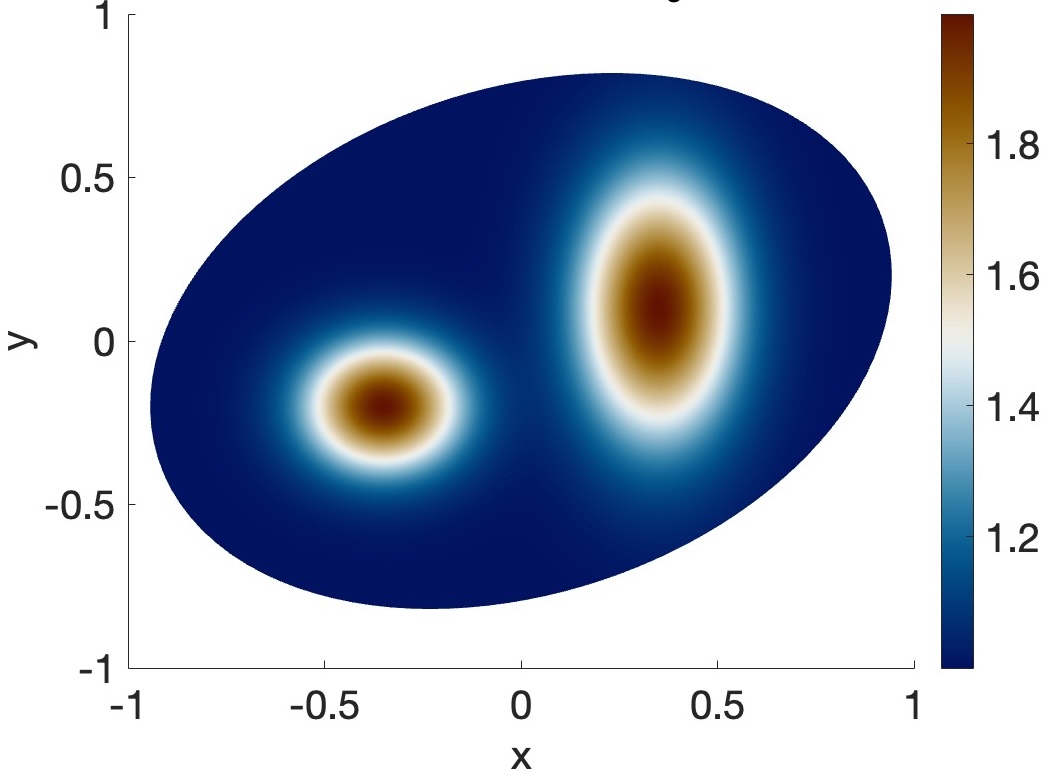}
\includegraphics[width=4.5cm,height=4cm]{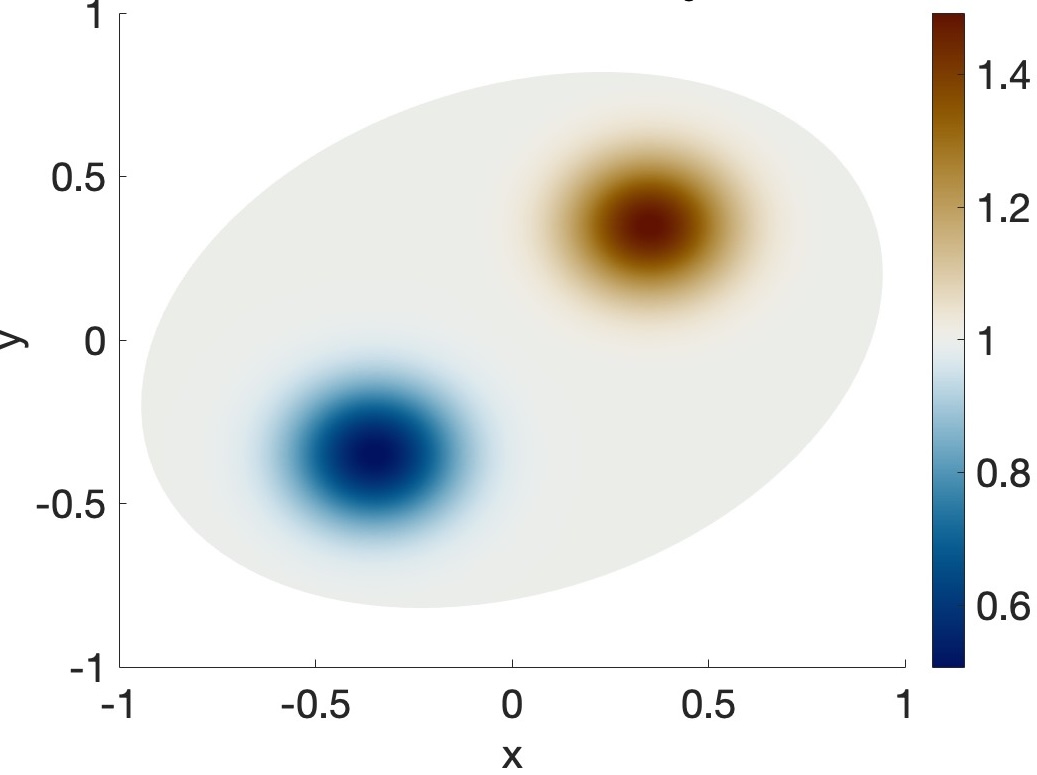}\\
\includegraphics[width=4.5cm,height=4cm]{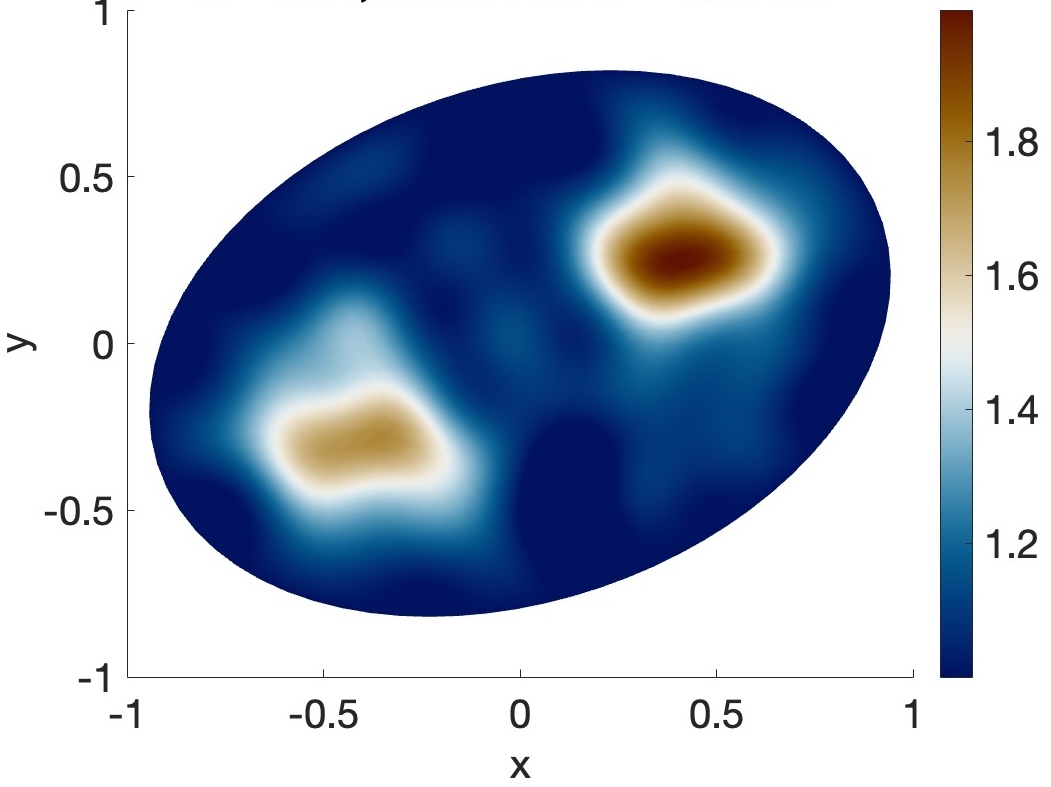}
\includegraphics[width=4.5cm,height=4cm]{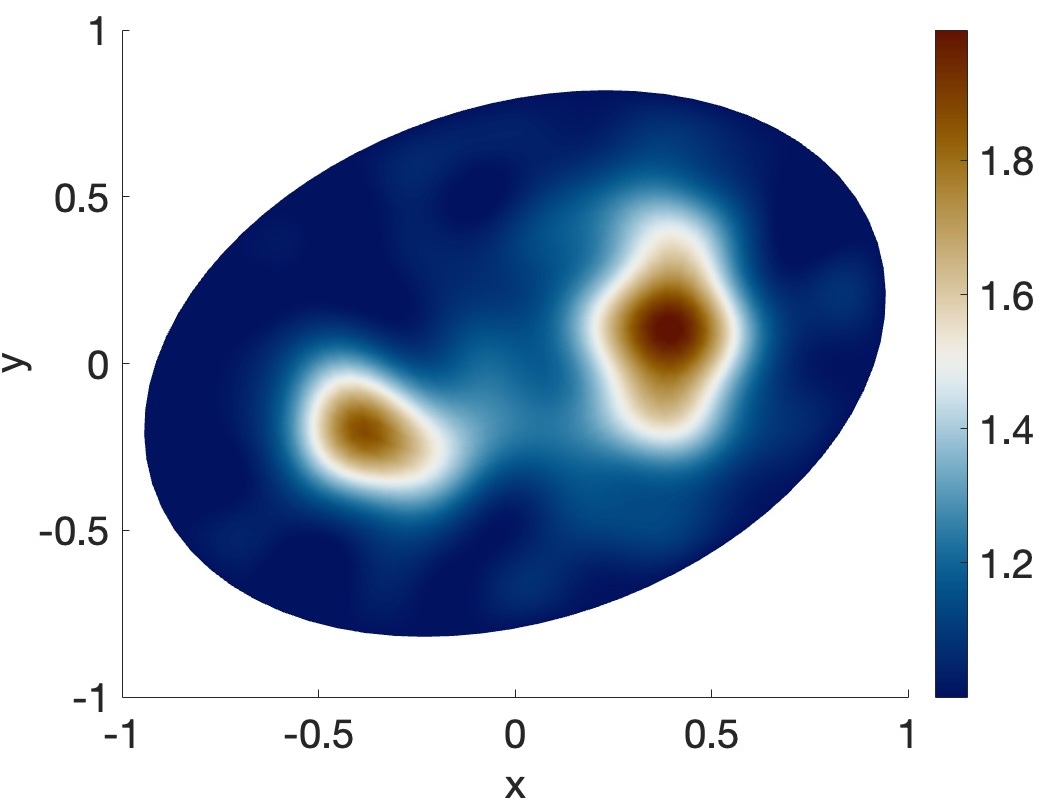}
\includegraphics[width=4.5cm,height=4cm]{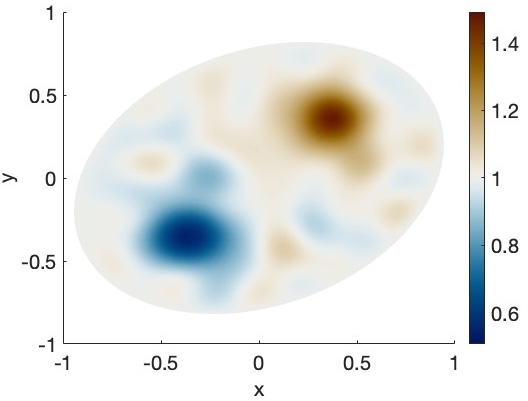}\\
\includegraphics[width=4.5cm,height=4cm]{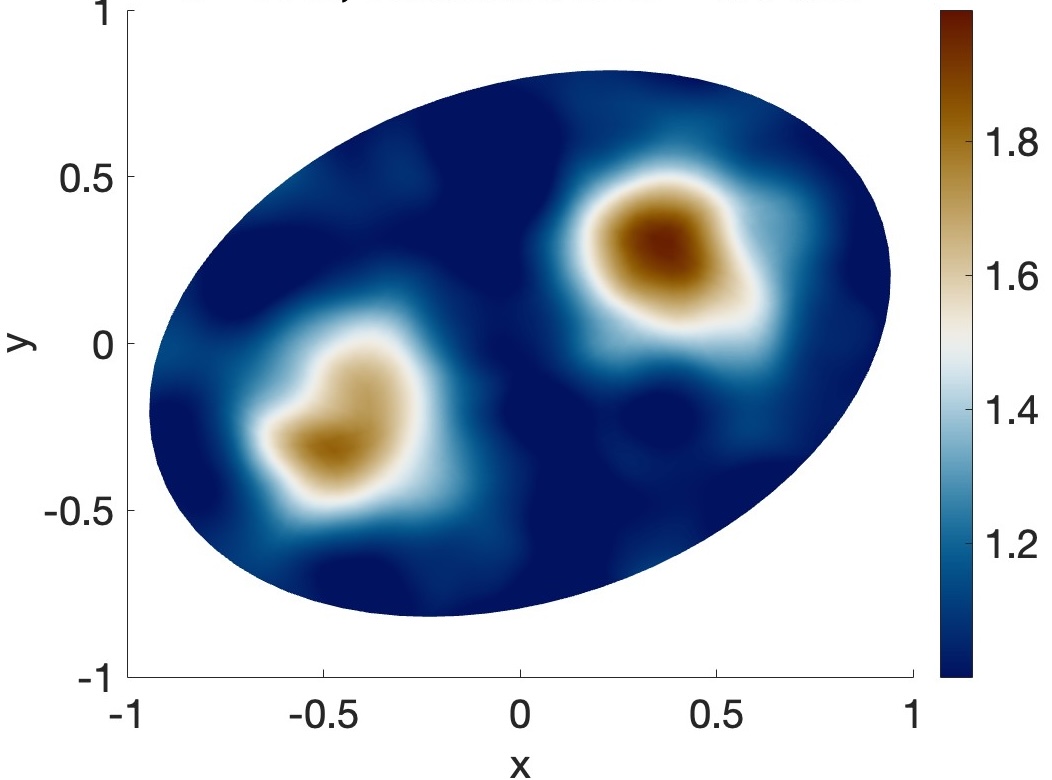}
\includegraphics[width=4.5cm,height=4cm]{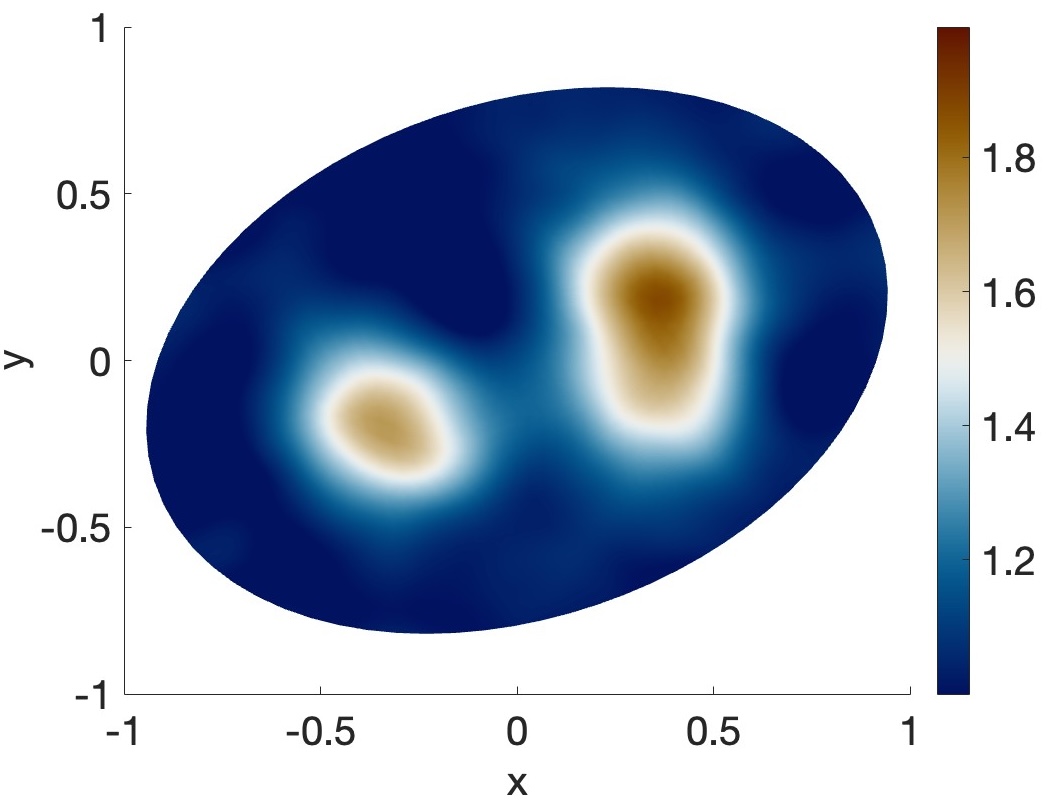}
\includegraphics[width=4.5cm,height=4cm]{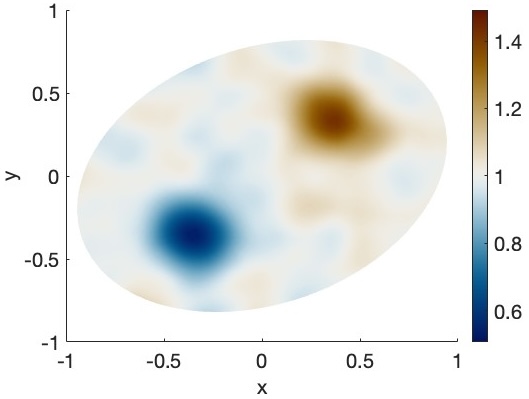}
\caption{Left %MDPI: Please change the hyphen (-) into minus sign ($-$, "U+2212"). e.g., "-1" should be "$-$1".
 column, top to bottom, respectively: the ground truth corresponding to $F_0^{(1)}$, and~the posterior mean estimates (computed via the pCN algorithm) for the truncated series and Matérn priors. The~step size is set to $\delta=0.00035,0.0001$ respectively, and~the acceptance ratios are $24.15\%$ and $26.51\%$. Central column: the ground truth for $F_0^{(2)}$ and the posterior mean estimates. For~pCN, $\delta=0.0001,0.00005$, acceptance ratios: $20.42\%$ and $22.73 \%$. Right column: the ground truth for $F_0^{(3)}$ and the posterior mean estimates. For~pCN, $\delta=0.0001,0.000035$, acceptance ratios: $25.61\%$, $29.85\%$. }
\label{Fig:NewEstimates}
\end{figure}
\unskip

\begin{table}[t]
\centering
\caption{Recovery %MDPI: %MDPI: 1. We removed the vertical lines. Please confirm. 2. Please check if the table has no table header.
 performances of the posterior mean estimate for different ground~truths.}
\label{Tab:NewResults}
%\centering
\renewcommand{\arraystretch}{1.75}
\begin{tabular}{cccc} 
\hline
 				 &   $F_0^{(1)}$  & $F_0^{(2)}$  & $F_0^{(3)}$  \\
\hline
 $\|\bar f_n - f_0\|_{L^2}$; series prior &  0.1197 &    0.0974  &  0.07234 \\
\hline
 $\|\bar f_n - f_0\|_{L^2}/\|f_0\|_{L^2}$; series prior &  7.22\% &  5.90\% &   4.71\% \\
\hline
 $\|\bar f_n - f_0\|_{L^2}$; Matérn prior & 0.1241 &    0.0941  &  0.08148  \\
\hline
 $\|\bar f_n - f_0\|_{L^2}/\|f_0\|_{L^2}$; Matérn prior &  7.49\% &  5.70\% &  5.31\% \\
\hline
\end{tabular}
\end{table}

We present further empirical investigations in which we consider the recovery of three additional true diffusivity functions, respectively specified (under the parametrisation \eqref{Eq:Param}) by
\begin{align*}
    F^{(1)}_0(x,y)&=  \log\left(1 + e^{-(4x-1.75)^2-(4y-1)^2}
    +e^{-(4x+1.75)^2-(4y+1)^2}\right);\\
    F^{(2)}_0(x,y)&= \log\left(1+ e^{-(5x-1.75)^2-(2.5y-0.25)^2}
	+e^{-(5x+1.75)^2-(5y+1)^2}\right);\\
	F^{(3)}_0(x,y)&=\log \left(1 + 0.5e^{-(5x-1.75)^2-(5y-1.75)^2)}
	-0.5e^{-(5x+1.75)^2-(5y+1.75)^2)}\right),
\end{align*}
for $(x,y)\in\Ocal$; see Figure~\ref{Fig:NewEstimates}, top row. For~each ground truth, we generate synthetic datasets $\{(Y_i,X_i)\}_{i=1}^n$, with~$n=2500$ as~described in Section~\ref{Sec:Results}, with~noise standard deviation $\sigma=0.0025$ (and signal-to-noise ratio respectively equal to 71.02, 73.26, and 63.88). Next, for~each set of observations, we implement posterior inference with truncated Gaussian series priors based on the Neumann--Laplacian eigenpairs and with the Matérn process prior, numerically computing the associated posterior mean estimates via the pCN algorithm. For~the series priors, the~regularity parameter in \eqref{Eq:DiscrPrior} is set to $\alpha=3$. For~the Matérn prior, the~covariance hyperparameters in \eqref{Eq:MatCov} are set to $\alpha = 3$ and $\ell = 0.25$. Across the three collections of experiments, the~pCN algorithm is iterated for 25,000 steps, with~a tuning of the step size to achieve a stabilisation of the acceptance rate between 20\% and 30\% after the burn-in (corresponding to the first 5000 samples). All the runs are initialised with cold starts. The~obtained results are visualised in Figure~\ref{Fig:NewEstimates} and summarised in Table~\ref{Tab:NewResults}. The~computation times are in line with those of the experiments presented in Section~\ref{Sec:Results}.

%
%
%
%
%

%%%%%%%%%%%%%%%%%%%%%%%%%%%%%%%%%%%%%%%%%%
\section[\appendixname~\thesection]{Proof of Theorem \ref{Theo:Main}}
\label{App:Proof}

We follow the proofs of Theorems 4--6 in~\cite{GN20}, recalling two key properties of the forward operator $G$, holding for all $\alpha>1+d/2$,
\begin{equation}
\label{Eq:LipEstim}
	\| G(\Phi\circ F_1) - G(\Phi\circ F_2)\|_{L^2}
	\lesssim (1 + \|F_1\|_{C^1}^4+\|F_2\|_{C^1}^4)\|F_1 - F_1\|_{(H^1)^*},
	\ \forall F_1,F_2\in H^\alpha_0(\Ocal),
\end{equation}
and
$$
	\sup_{F\in H^\alpha_0}\|G(\Phi\circ F)\|_{L^\infty}<\infty.
$$
\paragraph{{Step I:} posterior contraction rates in prediction risk. } We start with the derivation of the contraction rates in prediction risk \eqref{Eq:PredRiskCons}. Set $\epsilon_n := n^{-(\alpha+1)/(2\alpha+2+d)}$. By~an application of Theorem 13 and Lemmas 22 and 23 in~\cite{GN20}, it is enough to show that for some sufficiently large constant $c_1>0$, there exists a constant $c_2>0$ such that
\begin{equation}
\label{Eq:SmallBall}
	\Pi_n\left
	( F : \| G(\Phi\circ F) - G(\Phi\circ F_0)\|_{L^2} \le c_1 \epsilon_n
	\right)\ge e^{-c_2 n\epsilon_n^2},
\end{equation}
and further that there exist measurable sets $\Wcal_n\subseteq C^\beta(\Ocal)$ satisfying
\begin{equation}
\label{Eq:Sieves}
	\Pi_n(\Wcal_n^c)\le e^{-c_3 n\epsilon_n^2};
	\qquad \log N(\epsilon_n;\Wcal_n, d_{G})\lesssim n\epsilon_n^2,
\end{equation}
for some $c_3>0$ large enough, where $d_G(F_1,F_2):=\| G(\Phi\circ F_1) - G(\Phi\circ F_2)\|_{L^2}$, and $N(\epsilon_n;\Wcal_n, d_{G})$ is the minimal number of balls of radius $\epsilon_n$ in the metric $d_G$ needed to cover $\Acal_n$.

	The main difference compared to the proofs in~\cite{GN20} lies in the verification of the small ball probability estimate \eqref{Eq:SmallBall}. We proceed lower bounding this quantity by, for~sufficiently large $M>0$,
\begin{align*}
	\Pi_n\big(
	F : \| G(\Phi\circ F) &- G(\Phi\circ F_0)\|_{L^2} \le c_1 \epsilon_n, \|F - F_0\|_{C^1}\le M
	\big)\\
	&\ge \Pi_n\left( F : \|F-  F_0\|_{(H^1)^*} \le k_1 \epsilon_n, \|F - F_0\|_{C^1}\le M
	\right)
\end{align*}
for some $k_1>0$ (depending on $c_1$), having used \eqref{Eq:LipEstim}. Recalling the assumption  \eqref{Eq:ApproxSeq} on the approximating sequence $F_{0,n}\in \Hcal_{W_n}$ and using (twice) the triangle inequality, the~latter prior probability is greater than
\begin{align*}
	 \Pi_n\left( F : \|F-  F_{0,n}\|_{(H^1)^*} \le k_2 \epsilon_n, \|F - F_{0,n}\|_{C^1}\le M'
	\right)
	= \Pi_n\left(F : F- F_{0,n}\in B_1\cap B_2\right)
\end{align*}
for $k_2, M'>0$, having defined $B_1:=\{ F : \|F\|_{(H^1)^*} \le k_2 \epsilon_n\}$ and $B_2:=\{ F : \|F\|_{C_1} \le M'\}$. The~intersection $B_1\cap B_2$ defines a symmetric set in the ambient Banach space $C^1(\Ocal)$; hence, recalling that, by~linearity, the~RKHS of the scaled Gaussian prior $\Pi_n(\cdot)$ coincides with the RKHS $\Hcal_{W_n}$ of the base prior $\Pi_{W_n}(\cdot)$, with~the scaled RKHS norm being equal to $n^{d/(4\alpha+4+2d)}\|\cdot\|_{\Hcal_{W_n}}=\sqrt n \epsilon_n\|\cdot\|_{\Hcal_{W_n}}$ (e.g.,~Exercise 2.6.15 in~\cite{GN16}), an~application of Corollary 2.6.18 in~\cite{GN16} gives the further lower bound
$$
	e^{-n\epsilon_n^2\|F_{0,n}\|_{\Hcal_{W_n}}^2}\Pi_n(B_1\cap B_2)
	\ge e^{-k_3 n\epsilon_n^2} \Pi_n(B_1\cap B_2),
$$
for $k_3>0$ since $\|F_{0,n}\|_{\Hcal_{W_n}} = O(1)$ by assumption. Now, $B_1$ and $B_2$ are closed, convex, and symmetric subsets of the ambient space $C^1(\Ocal)$, and~therefore by the correlation inequality for Gaussian measures (cf.~Lemma A.2 in~\cite{giordano2022nonparametric}), 
$$
	\Pi_n (B_1\cap B_2)\ge \Pi_n (B_1)\Pi_n(B_2).
$$
Recalling the definition of the scaled Gaussian priors $\Pi_n(\cdot)$ in \eqref{Eq:Prior}, we have for $W\sim \Pi_{W_n}(\cdot)$,
$$
	\Pi_n(B_2) = \Pr (\|W\|_{C^1}\le M' \sqrt n \epsilon_n)
	\ge \Pr (\|W\|_{C^1}\le M'' ) >0
$$
for some $M''>0$ since $\sqrt n \epsilon_n >1$ for all $n$. Lastly, recalling the continuous embedding $\Hcal_{W_n}\subseteq H^\alpha_0(\Ocal)$ holding by assumption for some integer $\alpha>\beta+d/2$, $\beta\ge 1$, combining the metric entropy estimate
\begin{align*}
	\log &N(\eta; \{F\in\Hcal_{W_n} : \|F\|_{\Hcal_{W_n}}\le 1\}, \|\cdot\|_{(H^1)^*})\\
	&\le \log N(\eta; \{F\in H^\alpha_0(\Ocal) : \|F\|_{H^\alpha}\le k_4\}, \|\cdot\|_{(H^1)^*})
	\lesssim \eta^{-\frac{d}{\alpha+1}},
\end{align*}
cf.~Lemma 19 in~\cite{NvdGW20}, with~Theorem 1.2 in~\cite{LL99}, yields
\begin{align*}
	\Pi_n(B_1) 
	=\Pr(\|W\|_{(H^1)^*}\le k_2 \sqrt n \epsilon_n^2)
	\ge e^{-k_5 (\sqrt n \epsilon_n^2)^{-2\frac{d}{d+1}(2 - \frac{d}{\alpha+1})^{-1}}}
	= e^{-k_5 n\epsilon_n^2}.
\end{align*}
The obtained estimates thus jointly conclude the verification of \eqref{Eq:SmallBall} for a large enough constant $c_2$.

	Next, for~$K>0$, define the sieves
$$
	\Wcal_n :=\{ F: F = F_1 + F_2, 
	\ \|F_1\|_{(H^1)^*}\le K\epsilon_n, 
	\ \|F_2\|_{\Hcal_{W_n}}\le M, 
	\ \|F\|_{C^\beta}\le M\}.
$$
A direct application of Lemma 17 in~\cite{GN20} (with $\kappa=1$ in their notation) gives that for any positive $Q>0$, there exists sufficiently large $K$ such that $\Pi_n( \Wcal_n^c)\le e^{-Qn\epsilon_n^2}$. We then take $K$ that is large enough and argue as in the proof of Lemma 18 in~\cite{GN20} to deduce that $\log N(\epsilon_n;\Wcal_n,d_G)\lesssim n\epsilon_n^2$. This concludes the verification of \eqref{Eq:Sieves} and then also of the first claim \eqref{Eq:PredRiskCons} of Theorem \ref{Theo:Main} since by Theorem 13 in~\cite{GN20} we obtain that
\begin{equation}
\label{Eq:Interm}
	\Pi_n \Big(F\in\Wcal_n : \|G(\Phi\circ F)-G(\Phi\circ F)\|_{L^2}
	\le L\epsilon_n\Big|\{(Y_i,X_i)\}_{i=1}^n
	\Big) = 1 -O_{P_{f_0}^{(n)}}(e^{-D n\epsilon_n^2})
\end{equation}
as $n\to\infty$ for some $L, D>0$.

\paragraph{Step II: remaining claims.} Assume now that $\beta\ge2$, and~recall that $\alpha>\beta+d/2$ (with $\alpha,\beta\in\N$). Since $C^\beta(\Ocal)\subset H^\beta(\Ocal)$, Lemmas 23 and 29 in~\cite{NvdGW20} imply that for all $F\in \Wcal_n$, with~$\Wcal_n$ for the above sieve sets, we have $G(\Phi\circ F)\in H^{\beta+1}(\Ocal)$ and
$$
	\|G(\Phi\circ F)\|_{H^{\beta+1}}\lesssim 1 + \|\Phi\circ F\|_{H^\beta}^{\beta(\beta+1)}
	\lesssim 1,
$$
with a multiplicative constant independent of $F$. Similarly, since $F_0\in H^\alpha_0(\Ocal)\subset H^\beta(\Ocal)$, $G(\Phi\circ F_0)\in H^{\beta+1}(\Ocal)$ and $\|G(\Phi\circ F_0)\|_{H^{\beta+1}}\lesssim 1$. By~the standard interpolation inequality for Sobolev spaces, for~all $F\in\Wcal_n$, we then have
\begin{align*}
	\|G(\Phi\circ F) - G(\Phi\circ F_0) \|_{H^2}
	&\lesssim \|G(\Phi\circ F) - G(\Phi\circ F_0) \|_{L^2}^{\frac{\beta-1}{\beta+1}}
	\|G(\Phi\circ F) - G(\Phi\circ F_0) \|_{H^\beta}^\frac{2}{\beta+1}\\
	&\lesssim \|G(\Phi\circ F) - G(\Phi\circ F_0) \|_{L^2}^\frac{\beta-1}{\beta+1}.
\end{align*}
Combined with \eqref{Eq:Interm}, this implies that
$$
	\Pi_n \Big(F\in\Wcal_n : \|G(\Phi\circ F)-G(\Phi\circ F_0)\|_{H^2}
	\le L'\epsilon_n^{\frac{\beta-1}{\beta+1}}\Big|\{(Y_i,X_i)\}_{i=1}^n
	\Big) = 1 -O_{P_{f_0}^{(n)}}(e^{-D n\epsilon_n^2})
$$
for $L'>0$ as $n\to\infty$. The~second claim \eqref{Eq:Cons} of Theorem \ref{Theo:Main} is then verified using the following stability estimate for the forward operator $G$, 
$$
	\|f_1 - f_2\|_{L^2}\lesssim \|f\|_{C^1}\|G(f) - G(f_0)\|_{H^2}
$$
holding for all $f,f_0\in \Fcal_{\alpha,f_{\text{min}}}$ (the parameter space from \eqref{Eq:Param}) provided that $\alpha>2+d/2$ and $\inf_{x\in\Ocal}s(x)>0$, cf.~Lemma 24 in~\cite{NvdGW20}. Indeed, combined with the second to last display, this gives
$$
	\Pi_n \Big(F\in\Wcal_n : \|\Phi\circ F-\Phi\circ F_0\|_{L^2}
	\le L'\epsilon_n^{\frac{\beta-1}{\beta+1}}\Big|\{(Y_i,X_i)\}_{i=1}^n
	\Big) = 1 -O_{P_{f_0}^{(n)}}(e^{-D n\epsilon_n^2}).
$$
Based on the last display, the~proof of the last claim \eqref{Eq:ConvRate} now follows arguing exactly as for Theorem 6 in~\cite{GN20}.

%
%\begin{table}[H] 
%\caption{This is a table caption.\label{tab5}}
%%\newcolumntype{C}{>{\centering\arraybackslash}X}
%\begin{tabularx}{\textwidth}{CCC}
%\toprule
%\textbf{Title 1}	& \textbf{Title 2}	& \textbf{Title 3}\\
%\midrule
%Entry 1		& Data			& Data\\
%Entry 2		& Data			& Data\\
%\bottomrule
%\end{tabularx}
%\end{table}
%
%\section[\appendixname~\thesection]{}
%All appendix sections must be cited in the main text. In the appendices, Figures, Tables, etc. should be labeled, starting with ``A''---e.g., Figure A1, Figure A2, etc.

%%%%%%%%%%%%%%%%%%%%%%%%%%%%%%%%%%%%%%%%%%
%\isPreprints{}{% This command is only used for ``preprints''.
%\begin{adjustwidth}{-\extralength}{0cm}

%\reftitle{References}


\paragraph{Acknowledgement.} We are grateful to the Associate Editor and four anonymous referees for their helpful comments. This research has been partially supported by MUR, PRIN project 2022CLTYP4.

%\bibliography{References}

\bibliographystyle{acm}

\bibliography{PhDThesisReferences}

\begin{thebibliography}{999}

\bibitem[Engl et~al.(1996)Engl, Hanke, and Neubauer]{EHN96}
Engl, H.W.; Hanke, M.; Neubauer, A.
\newblock {\em Regularization of Inverse Problems}; Mathematics
  and Its Applications; Kluwer Academic Publishers Group: Dordrecht, The Netherlands, 1996; Volume 375, p. viii+321. %MDPI: Please check whether this is the correct page number.

\bibitem[Kaltenbacher et~al.(2008)Kaltenbacher, Neubauer, and Scherzer]{KNS08}
Kaltenbacher, B.; Neubauer, A.; Scherzer, O.
\newblock {\em Iterative Regularization Methods for Nonlinear Ill-Posed
  Problems}; {Radon Series on Computational and Applied
  Mathematics}; Walter de Gruyter GmbH \& Co. KG: Berlin, Germany, 2008; Volume~6, p. viii+194.
\newblock {\url{https://doi.org/10.1515/9783110208276}}.

\bibitem[Isakov(2017)]{I17}
Isakov, V.
\newblock {\em Inverse Problems for Partial Differential Equations}, third ed.; {Applied Mathematical Sciences}; Springer: Cham, Switzerland, 2017; Volume 127, p.~xv+406.
\newblock {\url{https://doi.org/10.1007/978-3-319-51658-5}}.

\bibitem[Kaipio and Somersalo(2004)]{KS04}
Kaipio, J.; Somersalo, E.
\newblock {\em Statistical and Computational Inverse Problems}; Number 160 in
  Applied Mathematical Sciences; Springer: New York, NY, USA, 2004.

\bibitem[Bissantz et~al.(2004)Bissantz, Hohage, and Munk]{BHM04}
Bissantz, N.; Hohage, T.; Munk, A.
\newblock Consistency and rates of convergence of nonlinear {T}ikhonov
  regularization with random noise.
\newblock {\em Inverse Probl.} {\bf 2004}, {\em 20},~1773--1789.

\bibitem[Hohage and Pricop(2008)]{HP08}
Hohage, T.; Pricop, M.
\newblock Nonlinear {T}ikhonov regularization in {H}ilbert scales for inverse
  boundary value problems with random noise.
\newblock {\em Inverse Probl. Imaging} {\bf 2008}, {\em 2},~271--290.
\newblock {\url{https://doi.org/10.3934/ipi.2008.2.271}}.

\bibitem[Benning and Burger(2018)]{BB18}
Benning, M.; Burger, M.
\newblock Modern regularization methods for inverse problems.
\newblock {\em Acta Numer.} {\bf 2018}, {\em 27},~1--111.
\newblock {\url{https://doi.org/10.1017/S0962492918000016}}.

\bibitem[Arridge et~al.(2019)Arridge, Maass, \"{O}ktem, and
  Sch\"{o}nlieb]{AMOS19}
Arridge, S.; Maass, P.; \"{O}ktem, O.; Sch\"{o}nlieb, C.B.
\newblock Solving inverse problems using data-driven models.
\newblock {\em Acta Numer.} {\bf 2019}, {\em 28},~1--174.
\newblock {\url{https://doi.org/10.1017/s0962492919000059}}.

\bibitem[Nickl(2023)]{N23}
Nickl, R.
\newblock {\em Bayesian Non-Linear Statistical Inverse Problems}; Zurich
  Lectures in Advanced Mathematics; EMS Press: Berlin, Germany, 2023; p. xi+159.
\newblock {\url{https://doi.org/10.4171/zlam/30}}.

\bibitem[Evans(2010)]{E10}
Evans, L.C.
\newblock {\em Partial Differential Equations}, 2nd ed.; {\em
  Graduate Studies in Mathematics}; American Mathematical Society: Providence,
  RI, USA, 2010; Volume~19, p. xxii+749.
\newblock {\url{https://doi.org/10.1090/gsm/019}}.

\bibitem[Yeh(1986)]{Y86}
Yeh, W.W.G.
\newblock Review of Parameter Identification Procedures in Groundwater
  Hydrology: The Inverse Problem.
\newblock {\em Water Resour. Res.} {\bf 1986}, {\em 22},~95--108.
%  \href{http://arxiv.org/abs/https://agupubs.onlinelibrary.wiley.com/doi/pdf/10.1029/WR022i002p00095}{{\normalfont
%  [https://agupubs.onlinelibrary.wiley.com/doi/pdf/10.1029/WR022i002p00095]}}.
\newblock {\url{https://doi.org/10.1029/WR022i002p00095}}.

\bibitem[Richter(1981)]{R81}
Richter, G.R.
\newblock An inverse problem for the steady state diffusion equation.
\newblock {\em SIAM J. Appl. Math.} {\bf 1981}, {\em 41},~210--221.
\newblock {\url{https://doi.org/10.1137/0141016}}.

\bibitem[Knowles(2001)]{K01}
Knowles, I.
\newblock Parameter identification for elliptic problems.
\newblock {\em J. Comput. Appl. Math.} {\bf 2001},
  {\em 131},~175--194.
\newblock
  {\url{https://doi.org/https://doi.org/10.1016/S0377-0427(00)00275-2}}.

\bibitem[Bonito et~al.(2017)Bonito, Cohen, DeVore, Petrova, and
  Welper]{BCDPW17}
Bonito, A.; Cohen, A.; DeVore, R.; Petrova, G.; Welper, G.
\newblock Diffusion coefficients estimation for elliptic partial differential
  equations.
\newblock {\em SIAM J. Math. Anal.} {\bf 2017}, {\em 49},~1570--1592.
\newblock {\url{https://doi.org/10.1137/16M1094476}}.

\bibitem[Dashti and Stuart(2011)]{DS11}
Dashti, M.; Stuart, A.M.
\newblock Uncertainty quantification and weak approximation of an elliptic
  inverse problem.
\newblock {\em SIAM J. Numer. Anal.} {\bf 2011}, {\em 49},~2524--2542.
\newblock {\url{https://doi.org/10.1137/100814664}}.

\bibitem[Cotter et~al.(2013)Cotter, Roberts, Stuart, and White]{CRSW13}
Cotter, S.; Roberts, G.; Stuart, A.; White, D.
\newblock {MCMC} Methods for Functions: Modifying Old Algorithms to Make Them
  Faster.
\newblock {\em Stat. Sci.} {\bf 2013}, {\em 28},~424--446.

\bibitem[Stuart(2010)]{S10}
Stuart, A.M.
\newblock Inverse problems: A {B}ayesian perspective.
\newblock {\em Acta Numer.} {\bf 2010}, {\em 19},~451--559.

\bibitem[Vollmer(2013)]{V13}
Vollmer, S.J.
\newblock Posterior consistency for {B}ayesian inverse problems through
  stability and regression results.
\newblock {\em Inverse Probl.} {\bf 2013}, {\em 29},~125011. %MDPI: We removed ``32'', please confirm.
\newblock {\url{https://doi.org/10.1088/0266-5611/29/12/125011}}.

\bibitem[Abraham and Nickl(2019)]{AN19}
Abraham, K.; Nickl, R.
\newblock On statistical {C}alder\'{o}n problems.
\newblock {\em Math. Stat. Learn.} {\bf 2019}, {\em 2},~165--216.

\bibitem[Giordano and Nickl(2020)]{GN20}
Giordano, M.; Nickl, R.
\newblock Consistency of {B}ayesian inference with {G}aussian process priors in
  an elliptic inverse problem.
\newblock {\em Inverse Probl.} {\bf 2020}, {\em 36},~85001--85036.

\bibitem[Monard et~al.(2021)Monard, Nickl, and Paternain]{MNP21a}
Monard, F.; Nickl, R.; Paternain, G.P.
\newblock Consistent inversion of noisy non-{A}belian {X}-ray transforms.
\newblock {\em Comm. Pure Appl. Math.} {\bf 2021}, {\em 74},~1045--1099.
\newblock {\url{https://doi.org/10.1002/cpa.21942}}.

\bibitem[Giordano(2021)]{G21}
Giordano, M.
\newblock Asymptotic Theory for Bayesian Nonparametric Inference in Statistical
  Models Arising from Partial Differential Equations.
\newblock {Doctoral Thesis, University of Cambridge, UK}, 2021. %MDPI: Please state the unversity and its location.
\newblock {\url{https://doi.org/10.17863/CAM.78120}}.

\bibitem[Agapiou and Wang(to appear)]{AW23}
Agapiou, S.; Wang, S.
\newblock Laplace priors and spatial inhomogeneity in Bayesian inverse
  problems.
\newblock {\em Bernoulli} \textbf{2024}, \emph{30}, 878--910.

\bibitem[Ghosal and van~der Vaart(2017)]{GvdV17}
Ghosal, S.; van~der Vaart, A.W.
\newblock {\em Fundamentals of Nonparametric Bayesian Inference}; Cambridge
  University Press: New York, NY, USA, 2017.

\bibitem[Gin\'e and Nickl(2016)]{GN16}
Gin\'e, E.; Nickl, R.
\newblock {\em Mathematical Foundations of Infinite-Dimensional Statistical
  Models}; Cambridge University Press: New York, NY, USA, 2016; p. xiv+690.
\newblock {\url{https://doi.org/10.1017/CBO9781107337862}}.

\bibitem[Nickl et~al.(2020)Nickl, van~de Geer, and Wang]{NvdGW20}
Nickl, R.; van~de Geer, S.; Wang, S.
\newblock Convergence rates for penalized least squares estimators in {PDE}
  constrained regression problems.
\newblock {\em SIAM/ASA J. Uncertain. Quantif.} {\bf 2020}, {\em 8},~374--413.
\newblock {\url{https://doi.org/10.1137/18M1236137}}.

\bibitem[Dashti and Stuart(2017)]{DS17}
Dashti, M.; Stuart, A.M.
\newblock The {B}ayesian approach to inverse problems. In {\em Handbook of
  Uncertainty Quantification}; Springer: Cham, Switzerland, 2017; 
 pp.~311--428.

\bibitem[Haroske and Triebel(2007)]{HT07}
Haroske, D.D.; Triebel, H.
\newblock {\em Distributions, Sobolev Spaces, Elliptic Equations}; EMS Press, Berlin, Germany %MDPI: Please add the location of the publisher (City and Country).
  2007.

\bibitem[Taylor(2011)]{T11}
Taylor, M.E.
\newblock {\em Partial Differential Equations I}; Springer: New York, NY, USA, 2011.

\bibitem[Hairer et~al.(2014)Hairer, Stuart, and Vollmer]{HSV14}
Hairer, M.; Stuart, A.; Vollmer, S.
\newblock Spectral Gaps for a {M}etropolis-{H}astings Algorithm in Infinite
  Dimensions.
\newblock {\em  Ann. Appl. Probab.} {\bf 2014}, {\em 24},~2455--2490.

\bibitem[Hall and Marron(1990)]{hall1990variance}
Hall, P.; Marron, J.
\newblock On variance estimation in nonparametric regression.
\newblock {\em Biometrika} {\bf 1990}, {\em 77},~415--419.

\bibitem[Wahba(1978)]{wahba1978improper}
Wahba, G.
\newblock Improper priors, spline smoothing and the problem of guarding against
  model errors in regression.
\newblock {\em J. R. Stat. Soc. Ser. B Stat. Methodol.} {\bf 1978}, {\em 40},~364--372.

\bibitem[Rice(1984)]{rice1984bandwidth}
Rice, J.
\newblock Bandwidth choice for nonparametric regression.
\newblock {\em Ann. Statist.}; 1984, 12,  %MDPI: Please add the name of the publisher and the location of it.
1215--1230.

\bibitem[Kejzlar et~al.(2021)Kejzlar, Son, Bhattacharya, and
  Maiti]{kejzlar2021fast}
Kejzlar, V.; Son, M.; Bhattacharya, S.; Maiti, T.
\newblock A fast and calibrated computer model emulator: An empirical Bayes
  approach.
\newblock {\em Stat. Comput.} {\bf 2021}, {\em 31},~49.

\bibitem[Knapik et~al.(2015)Knapik, Szab{\`o}, van~der Vaart, and van
  Zanten]{KSvdVvZ15}
Knapik, B.; Szab{\`o}, B.; van~der Vaart, A.W.; van Zanten, H.
\newblock Bayes procedures for adaptive inference in inverse problems for the
  white noise model.
\newblock {\em Probab. Theory Relat. Fields} {\textbf{2015}}, \emph{164}, 771--813.

\bibitem[Rousseau and Szabo(2017)]{roussszabo}
Rousseau, J.; Szabo, B.
\newblock {Asymptotic behaviour of the empirical Bayes posteriors associated to
  maximum marginal likelihood estimator}.
\newblock {\em  Ann. Stat.} {\bf 2017}, {\em 45},~833 -- 865.
\newblock {\url{https://doi.org/10.1214/16-AOS1469}}.

\bibitem[Teckentrup(2020)]{teckentrup2020convergence}
Teckentrup, A.L.
\newblock Convergence of Gaussian process regression with estimated
  hyper-parameters and applications in Bayesian inverse problems.
\newblock {\em SIAM/ASA J. Uncertain. Quantif.} {\bf 2020}, {\em
  8},~1310--1337.

\bibitem[Agapiou et~al.(2014)Agapiou, Bardsley, Papaspiliopoulos, and
  Stuart]{agapiou2014analysis}
Agapiou, S.; Bardsley, J.M.; Papaspiliopoulos, O.; Stuart, A.M.
\newblock Analysis of the Gibbs sampler for hierarchical inverse problems.
\newblock {\em SIAM/ASA J. Uncertain. Quantif.} {\bf 2014}, {\em 2},~511--544.

\bibitem[Giordano and Ray(2022)]{giordano2022nonparametric}
Giordano, M.; Ray, K.
\newblock Nonparametric Bayesian inference for reversible multidimensional
  diffusions.
\newblock {\em  Ann. Stat.} {\bf 2022}, {\em 50},~\mbox{2872--2898.}

\bibitem[Li and Linde(1999)]{LL99}
Li, W.V.; Linde, W.
\newblock Approximation, metric entropy and small ball estimates for {G}aussian
  measures.
\newblock {\em Ann. Probab.} {\bf 1999}, {\em 27},~1556--1578.
\newblock {\url{https://doi.org/10.1214/aop/1022677459}}.

\end{thebibliography}

\end{document}